\documentclass[11pt]{article}
\usepackage{amssymb, url, amsthm, amsmath}
\usepackage{amscd}
\usepackage[cp1251]{inputenc}
\usepackage[english,russian]{babel}
\newenvironment{proclaim}[1]{\bigskip \noindent \bf #1. \it}{\medskip}

\newcommand \N      {\mathbb N}
\newcommand \Z      {\mathbb Z}

\newcommand \ov     {\overline}
\newcommand \wh     {\widehat}
\newcommand \ssm    {\smallsetminus}
\newcommand \dd     {\partial}
\newcommand \la     {[}
\newcommand \ra     {]}
\newcommand \id     {\operatorname{id}}
\newcommand \si     {\sigma}
\newcommand \mi     {P_{\mu}}
\renewcommand \u    {\mathbf u}
\newcommand  \PP    {\mathcal P}
\newcommand  \qq    {\mathbf q}
\newcommand  \shl   {\operatorname{wing}}
\newcommand  \osn   {\operatorname{core}}

\newcommand  \mmBG  {1}
\newcommand  \mmRW  {2}
\newcommand  \mmRWL {3}
\newcommand  \mmPRIZ{4}
\newcommand  \mmF   {5}
\newcommand  \mmLF  {6}

\newcommand  \mmVYB {8}
\newcommand  \nJ    {10}

\urldef\vershik\url{vershik@pdmi.ras.ru} \urldef\malyutin\url{malyutin@pdmi.ras.ru}
\author {A.\,M.~Vershik, A.\,V.~Malyutin \thanks{%
St. Petersburg Branch, Steklov Mathematical Institute, Russian Academy of
Sciences, Fontanka 27, St. Petersburg 191023, Russia. E-mail: \vershik,
\malyutin. Partially supported by the grants RFFI 08-01-00379-а,
NSh-2460.2008.1112, CRDF RUM-2622ST-04. Keywords: braid group, random walk,
stable normal form, Poisson--Furstenberg boundary.}}

\date{26.04.07}
\title{\textbf{Boundaries of braid groups and\\ Markov--Ivanovsky normal form}}

\begin{document}
\maketitle
\renewcommand{\proofname}{Proof}
\renewcommand{\abstractname}{Abstract}
\renewcommand{\figurename}{Figure}
\renewcommand{\refname}{References}
\begin{abstract}
We describe random walk boundaries (in particular, the Poisson--Furstenberg, or
PF-boundary) for a vast family of groups in terms of the hyperbolic boundary
of a special free subgroup. We prove that almost all trajectories of the random
walk (with respect to an arbitrary nondegenerate measure on the group) converge
to points of that boundary. This implies the stability (in the sense
of~\cite{Ver}) of the so-called Markov--Ivanovsky normal form for braids.
\end{abstract}

\section*{Introduction}

\subsection{Random walks on groups and their boundaries}
In this paper we study boundaries of random walks on groups of a wide class,
the main examples of which are the Artin braid groups. Mainly we are interested
in the problem of {\it algebraic description of the boundary\/}; this problem
is still open during the last decades,~--- we construct and describe the
boundary of the random walk in terms of the group itself, or in other words in
terms of generators and relations. We do not prove that the boundaries which we
have found are {\it maximal\/}\footnote{Although for the braid groups this fact
follows from the comparison of our results with results of the
papers~\cite{KM,FM}.}. Nevertheless, the method we present joins Furstenberg's
approach to the boundaries with the idea of stable normal forms in groups, and
opens new possibilities for the description of boundaries.

The problem of the calculation of boundaries for countable groups (or more
generally, for graphs) consists of the following: consider the procedure of
sequential multiplication (for example, from the right side) of a current
element of the group $y_{n-1}\in G$, $n=1,2,\dots$, $y_0=id$, with a randomly
chosen (in accordance with a fixed probability distribution $\mu$ on the group)
element $x_n\in G$, $n=1,2,\dots$; how to describe those characteristics of the
growing product $y_n$ of random elements that are stabilized when $n$ tends to
infinity? For the braid groups this question is very natural in the case where
our fixed probability distribution is the uniform distribution on the set of
standard Artin's generators and their inverses: we choose randomly a
transposition of strands and multiply a current braid on this transposition;
what properties of the braid ``survived at infinity'' (such a properties we
will called ``stabilizing properties''). A priori it is not clear if such
nontrivial properties do exist (in other words, it is not clear if the boundary
is trivial or not).

There exist various approaches to the theory of boundaries of groups. First of
all, we can say about probabilistic boundaries, for example, {\it exit
boundaries\/} of Markov processes, or {\it Poisson boundaries\/} in the sense
of harmonic functions, etc. These boundaries are defined as measure spaces
without topology or canonical metric. For random walks, this boundary could be
defined as a quotient space over tail equivalence relation or other equivalence
relation. But such definition does not provide the direct description of the
boundary. Usually, when we consider a space which is a factor-space of the tail
equivalence relation, the main difficulty is to verify if this space is a
proper quotient of the boundary or the whole boundary itself. Up to now this
problem was solved for a very small class of groups. Technical question is: how
to define an appropriate metric or topology on the space which is supposed to
be boundary.

Another approach was given by H.\,Furstenberg~\cite{Fu1,Fu2,Fu3}. It starts
from the opposite point, namely from the study of topological or metric
$G$-spaces with quasi-invariant, stationary (with respect to the action of the
group $G$) measure. The group action on those spaces satisfies certain special
conditions: {\it $\mu$-proximality\/} or {\it mean-proximality\/} (see
Section~\mmRW\ for the definitions). Such spaces are called {\it
$\mu$-boundaries\/}, and probability boundaries (for example, the exit
boundary) if endowed with appropriate topology, satisfy these conditions. As
measure spaces, $\mu$-boundaries usually are quotient spaces of exit
boundaries. {\it A $\mu$-boundary that is isomorphic to the exit boundary is
called maximal\/}.\footnote{The exit boundary of a group can have several
distinct natural topological realizations. In other words, a group can have
several maximal $\mu$-boundaries. The braid group $B_n$ is a striking example.
From the one hand, the Thurston boundary of Teichm\"uller space of the sphere
with $n+1$ punctures (this boundary is homeomorphic to $(2n-5)$-sphere) is a
maximal $\mu$-boundary of $B_n$ (see~\cite{KM,FM} for the proof); from the
other hand, the exit boundary of $B_n$ can be realized as a circle; below, we
will describe the same boundary as the hyperbolic boundary of the free group of
rank $n-1$.} This approach seems to be very fruitful and became popular,
especially for classical groups and classical lattices. In the
survey~\cite{Ver} it was suggested to call the maximal $\mu$-boundary of the
random walk the {\it Poisson--Furstenberg boundary\/}; this term we use here.
For classical groups these boundaries usually can be described as boundaries of
some compactifications of the group. Thus, Furstenberg theory gives the
description of the boundaries as topological $G$-spaces with special measures.
This approach seems to be parallel to the theory of more subtle kind of
boundaries~--- {\it Martin boundaries\/}. But the problem of such approaches
lies in the absence of the direct link between $\mu$-boundary and the space of
trajectories of Markov process. In particular, no universal connection exists
between the convergence of trajectories and topology of the boundary.

Another one approach was proposed in the paper~\cite{Ver}~--- {\it it consists in
the construction of boundary as the limit space for normal forms of random
elements of the group}. More exactly, we try to find some kind of normal forms
for the group elements (for example, representing these elements as words over
an alphabet of generators of the group, or as geometrical configurations~---
``generalized words'', etc.) such that for almost all sequences of elements
there is a {\it convergence\/} of these normal forms to an infinite word or to
a limit configuration. Then we regard the space of infinite configurations as a
boundary of random walk. Examples: in the case of the free non-abelian group,
reduced words is a stable normal form, and the corresponding limit space is the
space of infinite reduced words; for the locally free groups, this is the space
of infinite heaps (see~\cite{BNV,Mal0}).
A more complicated example (from~\cite{Ver}) provides an interpretation of
boundaries of meta-abelian groups and wreath-products. Here we consider the
special type of configurations of paths on lattices as normal forms of the
elements of the meta-abelian groups, and stabilization means the convergence
for almost all sequences of configurations to an infinite configurations. Under
such approach, the connection between the convergence of trajectories of Markov
process and topology in the space of words underlies the notion of
stabilization. The main difficulty consists in the description of the action of
the group on the limit space, and in establishing the fact that the space
obtained is the boundary or even maximal boundary. The main result of this
paper shows that in the cases which we have considered, our approach could be
agreed with Furstenberg's approach and we can identify the space of limit
normal forms with PF-boundary.

Unfortunately, there are no universal tools to prove the boundary maximality;
this question remains open. There is the entropy criterion by
Kaimanovich~\cite{Kai}, which is a corollary of the global criteria of the
boundary triviality~\cite{KV,Der}, but it is difficult to apply it in concrete
situations, because for that we need an appropriate metric. Another plan, which
uses estimations of the stabilization rate of the normal forms, is still not
complete. In the case of mapping class groups, it was proved in~\cite{KM,FM}
(using the above-mentioned criteria) that the maximal boundary can be realized
as the Thurston boundary of Teichm\"uller space. Relatively cumbersome
combination, which we do not argue here, of those arguments and our results
allows us to claim that the boundary which we have found is maximal, but as it
was mentioned above, it is too indirect way, and we need further work in order
to make the algebraic theory of group boundaries self-consistent.

\subsection{Statement of results and a geometric illustration for the stability}
The results of the paper can be described as follows. First, for a countable
group with a normal free non-abelian subgroup it is proved that the ordinary
(hyperbolic) boundary of the subgroup is a $\mu$-boundary\footnote{In fact, in
the case of the braid group such a $\mu$-boundary is maximal (i.e., it is a
PF-boundary), providing that the measure $\mu$ has a finite first moment. But
this maximality follows, as mentioned above, by combining our results with the
results of~\cite{KM,FM}, which have a different nature. We think that our
approach leads to a direct proof of the boundary maximality. We conjecture that
such a boundary is maximal whenever the natural homomorphism of the group to
the group of automorphisms of the normal subgroup is injective.} of the group
itself (Theorem~1). Then, for a certain subclass (which is described below) of
groups with normal free subgroups it is shown that the projection of almost
every (a.e.) path into the normal free subgroup converges (with respect to the
hyperbolic compactification) to a point on the aforementioned hyperbolic
boundary (Theorem~2). Finally, implying these theorems to the braid group and
its {\it special\/} free subgroup, we prove that the Markov--Ivanovsky normal
form is stable (Theorem~3). With more details and with the exact wording, these
results are as follows.

\paragraph{The description of a $\mu$-boundary for the groups with a normal free subgroup.}
We recall that each automorphism of a {\it word-hyperbolic\/} group
$\mathcal{F}$ (take the free group as an example) can be uniquely extended to
a homeomorphism of its {\it hyperbolic boundary\/} $\dd \mathcal{F}$. In
particular, if $\mathcal{F}$ is a normal subgroup in a group $\mathcal{G}$,
then $\mathcal{G}$ acts on $\mathcal{F}$ by conjugation automorphisms, which
induces a natural continuous action of $\mathcal{G}$ on $\dd \mathcal{F}$.

\begin{proclaim}{Theorem~1}
Let $G$ be a countable group with a normal free non-abelian subgroup $F$. Let
$\mu$ be a nondegenerate measure on $G$. Then on the $G$-space $\dd F$ there
exists a unique $\mu$-stationary measure $\nu$, which is continuous. The pair
$(\dd F,\nu)$ is a $\mu$-boundary of $(G,\mu)$.
\end{proclaim}

In other words, Theorem~1 states that $\dd F$ is a {\it mean-proximal\/}
$G$-space, i.e., $\dd F$ is {\it $\mu$-proximal\/} for each nondegenerate
measure $\mu$ on $G$ (see Section~\mmRW\ or~\cite{Fu3} for the definitions).
The statement that $(\dd F,\nu)$ is a $\mu$-boundary of $(G,\mu)$ means by
definition that for a.e. path $\tau=\{\tau_i\}_{i\in\Z_+}$ of the right-hand
random $\mu$-walk on $G$ the sequence $\{\tau_i (\nu)\}_{i\in\Z_+}$ of
measures converges to a point measure $\delta_{w(\tau)}$ with $w(\tau)\in \dd
F$.

\paragraph{The convergence of paths for semidirect products.}
The problem of {\it pointwise convergence\/} for the mean-proximal $G$-space
$F\cup \dd F$ is the problem whether there exists a point $f\in F\cup \dd F$
such that for a.e. path $\tau=\{\tau_i\}_{i\in\Z_+}$ the sequence $\{\tau_i
(f)\}_{i\in\Z_+}$ converges to the above-defined point $w(\tau)\in \dd F$. In
the general case of Theorem~1, this problem is open. Nevertheless, we can
prove a pointwise convergence under certain additional assumptions on the
structure of our group.

\begin{proclaim}{Theorem~2}
Let $H$ be a countable group with a normal free non-abelian subgroup $F$.
Suppose that $H$ can be presented as a semidirect product of $F$ by a subgroup
$A$ {\rm(}i.e., $ H= F\rtimes A$\/{\rm)}. Suppose, moreover, that $F$ contains
a nontrivial element that is fixed under the action of the subgroup $A$. Then
for each nondegenerate measure $\mu$ on $H$ and for a.e. path
$\tau=\{\tau_i\}_{i\in\Z_+}=\{x_i\alpha_i\}_{i\in\Z_+}$ {\rm(}here, $x_i\in
F$, $\alpha_i\in A$\/{\rm)} of the random $\mu$-walk, the sequence
$\{x_i\}_{i\in\Z_+}$ of elements in $F$ converges {\rm(}in the hyperbolic
compactification $F\cup \dd F$\/{\rm)} to a point $w(\tau)$ on the boundary~$\dd
F$.
\end{proclaim}

\paragraph{The stability of the Markov--Ivanovsky normal form for the braid group.}
We can illustrate the above results by the example of the Artin braid group
$B_n$ (or by the slightly more general example of mapping class groups of
surfaces with nonempty boundary). It is well known that the pure braid group
$P_n$ (which has a finite index in $B_n$) is a semidirect product of its
normal free subgroup $F_{n-1}$ of rank $n-1$ by the pure $(n-1)$-braids
subgroup $P_{n-1}$. By applying Theorems~1 and~2 to this semidirect product,
we conclude that the hyperbolic boundary $\dd F_{n-1}$ is a $\mu$-boundary of
$P_n$ (in fact, it is maximal) and by projecting a path in $P_n$ to $F_{n-1}$
we almost surely (a.s.) obtain a convergent sequence (this sequence converges in the
hyperbolic metric to a point on $\dd F_{n-1}$). Since $P_n$ is a finite index
subgroup of $B_n$, we can extend these results to $B_n$. Now, we recall that
the definition of the Markov--Ivanovsky normal form \cite{Mar} based on the
decomposition $P_n=F_{n-1}\rtimes P_{n-1}$ (see Section~\mmBG\ for details).
Thus, in the case of the braid group, the above results convert into the
following one.

\begin{proclaim}{Theorem~3}
The Markov--Ivanovsky normal form in the Artin braid group is stable {\rm(}with
respect to the random walk with any nondegenerate distribution\/{\rm)}.
\end{proclaim}

Here, we give an example of a path of the random $\mu$-walk on the pure braid
group $P_4$, where $\mu$ is the uniform distribution on the set of generators
$s_{ji}$ and their inverses. (This set of generators is described in
Section~\mmBG\ below.) The $k$-th element of the path is denoted by $\gamma_k$.
On the right part of the list we write the Markov--Ivanovsky normal forms of
these elements. On the given part of the path the form stabilizes rapidly. As a
matter of fact, the lengths of the Markov--Ivanovsky normal forms for the
elements of a.e. path grow exponentially. Note that for a.e. path the
stabilization speed is also exponential.

\begin{equation*}
\begin{array}{ll}
\gamma_0=1;                               & \mathfrak{I}(\gamma_0)=1. \\
 \gamma_1=\gamma_0\cdot s_{31}^{-1};\quad & \mathfrak{I}(\gamma_1)=s_{31}^{-1}. \\
 \gamma_2=\gamma_1\cdot s_{41};           & \mathfrak{I}(\gamma_2)=s_{43}s_{41}s_{43}^{-1}s_{31}^{-1}. \\
 \gamma_3=\gamma_2\cdot s_{43}^{-1};      & \mathfrak{I}(\gamma_3)=s_{43}s_{41}s_{41}s_{43}^{-1}s_{41}^{-1}s_{43}^{-1}s_{31}^{-1}. \\
 \gamma_4=\gamma_3\cdot s_{32}^{-1};      & \mathfrak{I}(\gamma_4)=s_{43}s_{41}s_{41}s_{43}^{-1}s_{41}^{-1}s_{43}^{-1}s_{31}^{-1}s_{32}^{-1}. \\
 \gamma_5=\gamma_4\cdot s_{42};           & \mathfrak{I}(\gamma_5)=s_{43}s_{41}s_{41}s_{43}^{-1}s_{41}^{-1}s_{42}s_{43}^{-1}s_{31}^{-1}s_{32}^{-1}. \\
 \gamma_6=\gamma_5\cdot s_{21};           & \mathfrak{I}(\gamma_6)=s_{43}s_{41}s_{41}s_{43}^{-1}s_{41}^{-1}s_{42}s_{43}^{-1}s_{31}^{-1}s_{32}^{-1}s_{21}.\\
 \gamma_7=\gamma_6\cdot s_{32};           & \mathfrak{I}(\gamma_7)=s_{43}s_{41}s_{41}s_{43}^{-1}s_{41}^{-1}s_{42}s_{43}^{-1}s_{31}^{-1}s_{32}^{-1}s_{31}^{-1}s_{32}s_{31}s_{21}.\\
 \gamma_8=\gamma_7\cdot s_{41}^{-1};      & \mathfrak{I}(\gamma_8)=s_{43}s_{41}s_{41}s_{43}^{-1}s_{41}^{-1}s_{42}s_{41}^{-1}s_{42}^{-1}s_{41}^{-1}s_{42}s_{41}s_{43}^{-1}s_{31}^{-1}s_{32}^{-1}\cdots \\ 
 \gamma_9=\gamma_8\cdot s_{42}^{-1};      & \mathfrak{I}(\gamma_9)=s_{43}s_{41}s_{41}s_{43}^{-1}s_{41}^{-1}s_{42}s_{43}^{-1}s_{41}^{-1}s_{42}^{-1}s_{41}^{-1}s_{42}^{-1}s_{41}s_{42}s_{41}\cdots \\ 
\end{array}
\end{equation*}

In the definition of the Markov--Ivanovsky normal form, the decomposition
$P_m=F_{m-1}\rtimes P_{m-1}$ is employed step by step to the pure braid groups
with decreasing ranks. As a result, we obtain a decomposition of $P_n$ into a
product of $n-1$ free subgroups. Each braid $\gamma$ of $P_n$ is a unique
product of $n-1$ elements of these subgroups. The Markov--Ivanovsky normal
form $\mathfrak{I}(\gamma)$ of $\gamma$ is composed of $n-1$ parts that
correspond to these elements. (In the case of the group $P_4$, the form has
three parts. The first part uses the symbols, or the generators, $s_{43}$,
$s_{42}$, $s_{41}$, and their inverses. The second part includes $s_{32}$,
$s_{31}$, and inverses. The third part is a power of $s_{21}$. See the example
above.) For a.e. path of the random walk the length of each part of the form
(except the uttermost cyclic part) tends to infinity, so that the
stabilization of the form is the stabilization of its first part. It is an
interesting fact that other parts of the form (except the uttermost cyclic one
certainly) also stabilize for a.e. path, but in the limit all the stable
information about the path is reflected in the first part of the form.

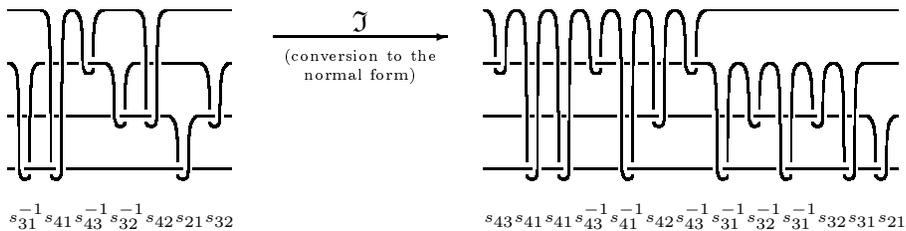
\begin{figure}[ht]
\centerline{\begin{picture}(336,115) \linethickness{0.7pt}
{\put(0,40){
\put(0,-20){$\scriptscriptstyle{s_{31}^{-1}s_{41}s_{43}^{-1}s_{32}^{-1}s_{42}s_{21}s_{32}}$}
\put(0,40) {\bezier29(0,0)(2,0)(3,-3)\put(3,-3){\bezier99(0,0)(1,-6)(1,-37.5)}
\put(4,-42.5){\bezier29(0,2)(0.5,-1.5)(2,-1.5)\bezier29(2,-1.5)(4,-1.5)(4,0)}
\put(8,-37){\bezier99(0,0)(0,28)(1,34)} \put(9,-3){\bezier19(0,0)(1,3)(3,3)}
\put(0,-40){\bezier9(0,0)(1,0)(2,0)\bezier19(6.5,0)(9,0)(12,0)}
\put(0,-20){\bezier9(0,0)(1,0)(1.5,0)\bezier9(5.7,0)(6,0)(6.3,0)\bezier9(10.2,0)(11,0)(12,0)}
\put(0,20){\bezier29(0,0)(6,0)(12,0)}}
\put(12,60) {\bezier29(0,0)(2,0)(3,-3)\put(3,-3){\bezier199(0,0)(1,-9)(1,-54)}
\put(4,-62.5){\bezier29(0,0)(0.5,-1.5)(2,-1.5)\bezier29(2,-1.5)(4,-1.5)(4,2)}
\put(8,-60.5){\bezier199(0,0)(0,48.5)(1,57.5)}
\put(9,-3){\bezier19(0,0)(1,3)(3,3)}
\put(0,-60){\bezier19(0,0)(5,0)(5.5,0)\bezier9(10,0)(11,0)(12,0)}
\put(0,-20){\bezier5(0,0)(1,0)(1.3,0)\bezier5(5.5,0)(6,0)(6.5,0)\bezier5(10.5,0)(11,0)(12,0)}
\put(0,-40){\bezier5(0,0)(1,0)(1.7,0)\bezier9(5.7,0)(6,0)(6.3,0)\bezier5(10.1,0)(11,0)(12,0)}
}
\put(24,60)
{\bezier29(0,0)(2,0)(3,-3)\put(3,-3){\bezier99(0,0)(1,-3)(1,-17.5)}
\put(4,-22.5){\bezier29(0,2)(0.5,-1.5)(2,-1.5)\bezier29(2,-1.5)(4,-1.5)(4,0)}
\put(8,-17){\bezier99(0,0)(0,11)(1,14)} \put(9,-3){\bezier19(0,0)(1,3)(3,3)}
\put(0,-20){\bezier9(0,0)(1,0)(2,0)\bezier19(6.5,0)(9,0)(12,0)}
\put(0,-40){\bezier29(0,0)(6,0)(12,0)}\put(0,-60){\bezier29(0,0)(6,0)(12,0)}}
\put(36,40)
{\bezier29(0,0)(2,0)(3,-3)\put(3,-3){\bezier99(0,0)(1,-3)(1,-17.5)}
\put(4,-22.5){\bezier29(0,2)(0.5,-1.5)(2,-1.5)\bezier29(2,-1.5)(4,-1.5)(4,0)}
\put(8,-17){\bezier99(0,0)(0,11)(1,14)} \put(9,-3){\bezier19(0,0)(1,3)(3,3)}
\put(0,-20){\bezier9(0,0)(1,0)(2,0)\bezier19(6.5,0)(9,0)(12,0)}
\put(0,20){\bezier29(0,0)(6,0)(12,0)}\put(0,-40){\bezier29(0,0)(6,0)(12,0)}}
\put(48,60) {\bezier29(0,0)(2,0)(3,-3)\put(3,-3){\bezier99(0,0)(1,-6)(1,-34)}
\put(4,-42.5){\bezier29(0,0)(0.5,-1.5)(2,-1.5)\bezier29(2,-1.5)(4,-1.5)(4,2)}
\put(8,-40.5){\bezier99(0,0)(0,31.5)(1,37.5)}
\put(9,-3){\bezier19(0,0)(1,3)(3,3)}
\put(0,-40){\bezier19(0,0)(5,0)(5.5,0)\bezier9(10,0)(11,0)(12,0)}
\put(0,-20){\bezier9(0,0)(1,0)(1.5,0)\bezier9(5.7,0)(6,0)(6.3,0)\bezier9(10.2,0)(11,0)(12,0)}
\put(0,-60){\bezier29(0,0)(6,0)(12,0)}}
\put(60,20) {\bezier29(0,0)(2,0)(3,-3)\put(3,-3){\bezier99(0,0)(1,-3)(1,-14)}
\put(4,-22.5){\bezier29(0,0)(0.5,-1.5)(2,-1.5)\bezier29(2,-1.5)(4,-1.5)(4,2)}
\put(8,-20.5){\bezier99(0,0)(0,14.5)(1,17.5)}
\put(9,-3){\bezier19(0,0)(1,3)(3,3)}
\put(0,-20){\bezier19(0,0)(5,0)(5.5,0)\bezier9(10,0)(11,0)(12,0)}
\put(0,20){\bezier29(0,0)(6,0)(12,0)}\put(0,40){\bezier29(0,0)(6,0)(12,0)}}
\put(72,40) {\bezier29(0,0)(2,0)(3,-3)\put(3,-3){\bezier99(0,0)(1,-3)(1,-14)}
\put(4,-22.5){\bezier29(0,0)(0.5,-1.5)(2,-1.5)\bezier29(2,-1.5)(4,-1.5)(4,2)}
\put(8,-20.5){\bezier99(0,0)(0,14.5)(1,17.5)}
\put(9,-3){\bezier19(0,0)(1,3)(3,3)}
\put(0,-20){\bezier19(0,0)(5,0)(5.5,0)\bezier9(10,0)(11,0)(12,0)}
\put(0,-40){\bezier29(0,0)(6,0)(12,0)}\put(0,20){\bezier29(0,0)(6,0)(12,0)}}
}
\put(100,90){\vector(1,0){66}} \put(100,90){\makebox(66,12){$\mathfrak{I}$}}
\put(100,90){\makebox(66,-15){\tiny (conversion to the}}
\put(100,90){\makebox(66,-30){\tiny normal form)}}
\put(180,40){
\put(0,-20){$\scriptscriptstyle{s_{43}s_{41}s_{41}s_{43}^{-1}s_{41}^{-1}s_{42}s_{43}^{-1}s_{31}^{-1}s_{32}^{-1}s_{31}^{-1}s_{32}s_{31}s_{21}}$}
\put(0,60) {\bezier29(0,0)(2,0)(3,-3)\put(3,-3){\bezier99(0,0)(1,-3)(1,-14)}
\put(4,-22.5){\bezier29(0,0)(0.5,-1.5)(2,-1.5)\bezier29(2,-1.5)(4,-1.5)(4,2)}
\put(8,-20.5){\bezier99(0,0)(0,14.5)(1,17.5)}
\put(9,-3){\bezier19(0,0)(1,3)(3,3)}
\put(0,-20){\bezier19(0,0)(5,0)(5.5,0)\bezier9(10,0)(11,0)(12,0)}
\put(0,-40){\bezier29(0,0)(6,0)(12,0)}\put(0,-60){\bezier29(0,0)(6,0)(12,0)}}
\put(12,60) {\bezier29(0,0)(2,0)(3,-3)\put(3,-3){\bezier199(0,0)(1,-9)(1,-54)}
\put(4,-62.5){\bezier29(0,0)(0.5,-1.5)(2,-1.5)\bezier29(2,-1.5)(4,-1.5)(4,2)}
\put(8,-60.5){\bezier199(0,0)(0,48.5)(1,57.5)}
\put(9,-3){\bezier19(0,0)(1,3)(3,3)}
\put(0,-60){\bezier19(0,0)(5,0)(5.5,0)\bezier9(10,0)(11,0)(12,0)}
\put(0,-20){\bezier5(0,0)(1,0)(1.3,0)\bezier5(5.5,0)(6,0)(6.5,0)\bezier5(10.5,0)(11,0)(12,0)}
\put(0,-40){\bezier5(0,0)(1,0)(1.7,0)\bezier9(5.7,0)(6,0)(6.3,0)\bezier5(10.1,0)(11,0)(12,0)}
}
\put(24,60) {\bezier29(0,0)(2,0)(3,-3)\put(3,-3){\bezier199(0,0)(1,-9)(1,-54)}
\put(4,-62.5){\bezier29(0,0)(0.5,-1.5)(2,-1.5)\bezier29(2,-1.5)(4,-1.5)(4,2)}
\put(8,-60.5){\bezier199(0,0)(0,48.5)(1,57.5)}
\put(9,-3){\bezier19(0,0)(1,3)(3,3)}
\put(0,-60){\bezier19(0,0)(5,0)(5.5,0)\bezier9(10,0)(11,0)(12,0)}
\put(0,-20){\bezier5(0,0)(1,0)(1.3,0)\bezier5(5.5,0)(6,0)(6.5,0)\bezier5(10.5,0)(11,0)(12,0)}
\put(0,-40){\bezier5(0,0)(1,0)(1.7,0)\bezier9(5.7,0)(6,0)(6.3,0)\bezier5(10.1,0)(11,0)(12,0)}
}
\put(36,60)
{\bezier29(0,0)(2,0)(3,-3)\put(3,-3){\bezier99(0,0)(1,-3)(1,-17.5)}
\put(4,-22.5){\bezier29(0,2)(0.5,-1.5)(2,-1.5)\bezier29(2,-1.5)(4,-1.5)(4,0)}
\put(8,-17){\bezier99(0,0)(0,11)(1,14)} \put(9,-3){\bezier19(0,0)(1,3)(3,3)}
\put(0,-20){\bezier9(0,0)(1,0)(2,0)\bezier19(6.5,0)(9,0)(12,0)}
\put(0,-40){\bezier29(0,0)(6,0)(12,0)}\put(0,-60){\bezier29(0,0)(6,0)(12,0)}}
\put(48,60)
{\bezier29(0,0)(2,0)(3,-3)\put(3,-3){\bezier199(0,0)(1,-9)(1,-57.5)}
\put(4,-62.5){\bezier29(0,2)(0.5,-1.5)(2,-1.5)\bezier29(2,-1.5)(4,-1.5)(4,0)}
\put(8,-57){\bezier199(0,0)(0,45)(1,54)} \put(9,-3){\bezier19(0,0)(1,3)(3,3)}
\put(0,-60){\bezier9(0,0)(1,0)(2,0)\bezier19(6.5,0)(9,0)(12,0)}
\put(0,-20){\bezier5(0,0)(1,0)(1.3,0)\bezier5(5.5,0)(6,0)(6.5,0)\bezier5(10.5,0)(11,0)(12,0)}
\put(0,-40){\bezier5(0,0)(1,0)(1.7,0)\bezier9(5.7,0)(6,0)(6.3,0)\bezier5(10.1,0)(11,0)(12,0)}
}
\put(60,60) {\bezier29(0,0)(2,0)(3,-3)\put(3,-3){\bezier99(0,0)(1,-6)(1,-34)}
\put(4,-42.5){\bezier29(0,0)(0.5,-1.5)(2,-1.5)\bezier29(2,-1.5)(4,-1.5)(4,2)}
\put(8,-40.5){\bezier99(0,0)(0,31.5)(1,37.5)}
\put(9,-3){\bezier19(0,0)(1,3)(3,3)}
\put(0,-40){\bezier19(0,0)(5,0)(5.5,0)\bezier9(10,0)(11,0)(12,0)}
\put(0,-20){\bezier9(0,0)(1,0)(1.5,0)\bezier9(5.7,0)(6,0)(6.3,0)\bezier9(10.2,0)(11,0)(12,0)}
\put(0,-60){\bezier29(0,0)(6,0)(12,0)}}
\put(72,60)
{\bezier29(0,0)(2,0)(3,-3)\put(3,-3){\bezier99(0,0)(1,-3)(1,-17.5)}
\put(4,-22.5){\bezier29(0,2)(0.5,-1.5)(2,-1.5)\bezier29(2,-1.5)(4,-1.5)(4,0)}
\put(8,-17){\bezier99(0,0)(0,11)(1,14)} \put(9,-3){\bezier19(0,0)(1,3)(3,3)}
\put(0,-20){\bezier9(0,0)(1,0)(2,0)\bezier19(6.5,0)(9,0)(12,0)}
\put(0,-40){\bezier29(0,0)(6,0)(12,0)}\put(0,-60){\bezier29(0,0)(6,0)(12,0)}}
\put(84,40)
{\bezier29(0,0)(2,0)(3,-3)\put(3,-3){\bezier99(0,0)(1,-6)(1,-37.5)}
\put(4,-42.5){\bezier29(0,2)(0.5,-1.5)(2,-1.5)\bezier29(2,-1.5)(4,-1.5)(4,0)}
\put(8,-37){\bezier99(0,0)(0,28)(1,34)} \put(9,-3){\bezier19(0,0)(1,3)(3,3)}
\put(0,-40){\bezier9(0,0)(1,0)(2,0)\bezier19(6.5,0)(9,0)(12,0)}
\put(0,-20){\bezier9(0,0)(1,0)(1.5,0)\bezier9(5.7,0)(6,0)(6.3,0)\bezier9(10.2,0)(11,0)(12,0)}
\put(0,20){\bezier29(0,0)(6,0)(12,0)}}
\put(96,40)
{\bezier29(0,0)(2,0)(3,-3)\put(3,-3){\bezier99(0,0)(1,-3)(1,-17.5)}
\put(4,-22.5){\bezier29(0,2)(0.5,-1.5)(2,-1.5)\bezier29(2,-1.5)(4,-1.5)(4,0)}
\put(8,-17){\bezier99(0,0)(0,11)(1,14)} \put(9,-3){\bezier19(0,0)(1,3)(3,3)}
\put(0,-20){\bezier9(0,0)(1,0)(2,0)\bezier19(6.5,0)(9,0)(12,0)}
\put(0,20){\bezier29(0,0)(6,0)(12,0)}\put(0,-40){\bezier29(0,0)(6,0)(12,0)}}
\put(108,40)
{\bezier29(0,0)(2,0)(3,-3)\put(3,-3){\bezier99(0,0)(1,-6)(1,-37.5)}
\put(4,-42.5){\bezier29(0,2)(0.5,-1.5)(2,-1.5)\bezier29(2,-1.5)(4,-1.5)(4,0)}
\put(8,-37){\bezier99(0,0)(0,28)(1,34)} \put(9,-3){\bezier19(0,0)(1,3)(3,3)}
\put(0,-40){\bezier9(0,0)(1,0)(2,0)\bezier19(6.5,0)(9,0)(12,0)}
\put(0,-20){\bezier9(0,0)(1,0)(1.5,0)\bezier9(5.7,0)(6,0)(6.3,0)\bezier9(10.2,0)(11,0)(12,0)}
\put(0,20){\bezier29(0,0)(6,0)(12,0)}}
\put(120,40) {\bezier29(0,0)(2,0)(3,-3)\put(3,-3){\bezier99(0,0)(1,-3)(1,-14)}
\put(4,-22.5){\bezier29(0,0)(0.5,-1.5)(2,-1.5)\bezier29(2,-1.5)(4,-1.5)(4,2)}
\put(8,-20.5){\bezier99(0,0)(0,14.5)(1,17.5)}
\put(9,-3){\bezier19(0,0)(1,3)(3,3)}
\put(0,-20){\bezier19(0,0)(5,0)(5.5,0)\bezier9(10,0)(11,0)(12,0)}
\put(0,-40){\bezier29(0,0)(6,0)(12,0)}\put(0,20){\bezier29(0,0)(6,0)(12,0)}}
\put(132,40) {\bezier29(0,0)(2,0)(3,-3)\put(3,-3){\bezier99(0,0)(1,-6)(1,-34)}
\put(4,-42.5){\bezier29(0,0)(0.5,-1.5)(2,-1.5)\bezier29(2,-1.5)(4,-1.5)(4,2)}
\put(8,-40.5){\bezier99(0,0)(0,31.5)(1,37.5)}
\put(9,-3){\bezier19(0,0)(1,3)(3,3)}
\put(0,-40){\bezier19(0,0)(5,0)(5.5,0)\bezier9(10,0)(11,0)(12,0)}
\put(0,-20){\bezier9(0,0)(1,0)(1.5,0)\bezier9(5.7,0)(6,0)(6.3,0)\bezier9(10.2,0)(11,0)(12,0)}
\put(0,20){\bezier29(0,0)(6,0)(12,0)}}
\put(144,20) {\bezier29(0,0)(2,0)(3,-3)\put(3,-3){\bezier99(0,0)(1,-3)(1,-14)}
\put(4,-22.5){\bezier29(0,0)(0.5,-1.5)(2,-1.5)\bezier29(2,-1.5)(4,-1.5)(4,2)}
\put(8,-20.5){\bezier99(0,0)(0,14.5)(1,17.5)}
\put(9,-3){\bezier19(0,0)(1,3)(3,3)}
\put(0,-20){\bezier19(0,0)(5,0)(5.5,0)\bezier9(10,0)(11,0)(12,0)}
\put(0,20){\bezier29(0,0)(6,0)(12,0)}\put(0,40){\bezier29(0,0)(6,0)(12,0)}}
}}
\end{picture}}
\caption{Geometric presentations for the pure $4$-braid
$s_{31}^{-1}s_{41}s_{43}^{-1}s_{32}^{-1}s_{42}s_{21}s_{32}$
and its Markov--Ivanovsky normal form.} \label{geom-intro}
\end{figure}

If we consider the ordinary geometric interpretation of braids, where a braid
is represented by a collections of intertwining strings, then the
Markov--Ivanovsky normal form of a braid is represented by a strings
interweaving of the following shape. First, the outermost string wind round
other strings, which are motionless at this moment (this is the first part of
the form). Then, only the second string moves (the second part of the form),
and so on. In the last (cyclic) part of the form, only two strings are
twisting. See Figure~\ref{geom-intro}; the generator $s_{ji}$ corresponds on
the geometric level to the one full twist between $i$-th and $j$-th strings.

Figure~\ref{geom-path} shows the geometric shapes for the normal forms of
elements in the above example. In this geometric interpretation, the
stabilization of the Markov--Ivanovsky normal form is observed as the
stabilization of the first string's position.

\begin{figure}[ht]
\unitlength=1.3pt \centerline{\begin{picture}(300,350)
\linethickness{0,5pt}\put(150,0){
\put(-43,315) {
\put(-93,16){$\gamma_0=1;$}\put(-30,16){$\mathfrak{I}(\gamma_0)=$}
\put(0,0){\put(1,16){$1.$}\multiput(0,0)(0,4){4}{\bezier29(0,0)(6,0)(12,0)}}
}
\put(-43,280) {
\put(-93,16){$\gamma_1=\gamma_0\cdot
s_{31}^{-1};$}\put(-30,16){$\mathfrak{I}(\gamma_1)=$}
\put(0,0){\put(1,16){$s_{31}^{-1}.$}\put(0,8){\put(0,0){\bezier29(0,0)(4,0)(4,-6.3)}\put(4,-8){\bezier19(0,1.7)(0,-2.5)(3.5,-1.5)}\put(8,0){\bezier29(0,-6.3)(0,0)(4,0)}}\put(0,0){\put(0,0){\bezier19(0,0)(1,0)(2.1,0)}\put(12,0){\bezier19(-5.7,0)(-3,0)(0,0)}}\put(0,4){\put(0,0){\bezier9(0,0)(1,0)(1.7,0)}\put(6,0){\bezier9(-0.3,0)(0,0)(0.3,0)}\put(12,0){\bezier9(-1.7,0)(-1,0)(0,0)}}\put(0,12){\bezier29(0,0)(6,0)(12,0)}}
}
\put(-43,245) {
\put(-93,16){$\gamma_2=\gamma_1\cdot
s_{41};$}\put(-30,16){$\mathfrak{I}(\gamma_2)=$}
\put(0,0){\put(1,16){$s_{43}$}\put(0,12){\put(0,0){\bezier19(0,0)(4,0)(4,-2.3)}\put(4,-4){\bezier19(0.5,-1.5)(4,-2.5)(4,1.7)}\put(8,0){\bezier19(0,-2.3)(0,0)(4,0)}}\put(0,8){\put(0,0){\bezier19(0,0)(3,0)(5.7,0)}\put(12,0){\bezier9(-2.1,0)(-1,0)(0,0)}}\put(0,0){\bezier29(0,0)(6,0)(12,0)}\put(0,4){\bezier29(0,0)(6,0)(12,0)}}
\put(12,0){\put(1,16){$s_{41}$}\put(0,12){\put(0,0){\bezier39(0,0)(4,0)(4,-10.3)}\put(4,-12){\bezier19(0.5,-1.5)(4,-2.5)(4,1.7)}\put(8,0){\bezier39(0,-10.3)(0,0)(4,0)}}\put(0,0){\put(0,0){\bezier19(0,0)(3,0)(5.7,0)}\put(12,0){\bezier19(-2.1,0)(-1,0)(0,0)}}\put(0,4){\put(0,0){\bezier9(0,0)(1,0)(2,0)}\put(6,0){\bezier9(-0.2,0)(0,0)(0.2,0)}\put(12,0){\bezier9(-2,0)(-1,0)(0,0)}}\put(0,8){\put(0,0){\bezier9(0,0)(1,0)(1.5,0)}\put(6,0){\bezier9(-0.7,0)(0,0)(0.7,0)}\put(12,0){\bezier9(-1.5,0)(-1,0)(0,0)}}}
\put(24,0){\put(1,16){$s_{43}^{-1}$}\put(0,12){\put(0,0){\bezier19(0,0)(4,0)(4,-2.3)}\put(4,-4){\bezier19(0,1.7)(0,-2.5)(3.5,-1.5)}\put(8,0){\bezier19(0,-2.3)(0,0)(4,0)}}\put(0,8){\put(0,0){\bezier19(0,0)(1,0)(2.1,0)}\put(12,0){\bezier19(-5.7,0)(-3,0)(0,0)}}\put(0,0){\bezier29(0,0)(6,0)(12,0)}\put(0,4){\bezier29(0,0)(6,0)(12,0)}}
\put(36,0){\put(1,16){$s_{31}^{-1}.$}\put(0,8){\put(0,0){\bezier29(0,0)(4,0)(4,-6.3)}\put(4,-8){\bezier19(0,1.7)(0,-2.5)(3.5,-1.5)}\put(8,0){\bezier29(0,-6.3)(0,0)(4,0)}}\put(0,0){\put(0,0){\bezier19(0,0)(1,0)(2.1,0)}\put(12,0){\bezier19(-5.7,0)(-3,0)(0,0)}}\put(0,4){\put(0,0){\bezier9(0,0)(1,0)(1.7,0)}\put(6,0){\bezier9(-0.3,0)(0,0)(0.3,0)}\put(12,0){\bezier9(-1.7,0)(-1,0)(0,0)}}\put(0,12){\bezier29(0,0)(6,0)(12,0)}}
}
\put(-43,210) {
\put(-93,16){$\gamma_3=\gamma_2\cdot
s_{43}^{-1};$}\put(-30,16){$\mathfrak{I}(\gamma_3)=$}
\put(0,0){\put(1,16){$s_{43}$}\put(0,12){\put(0,0){\bezier19(0,0)(4,0)(4,-2.3)}\put(4,-4){\bezier19(0.5,-1.5)(4,-2.5)(4,1.7)}\put(8,0){\bezier19(0,-2.3)(0,0)(4,0)}}\put(0,8){\put(0,0){\bezier19(0,0)(3,0)(5.7,0)}\put(12,0){\bezier9(-2.1,0)(-1,0)(0,0)}}\put(0,0){\bezier29(0,0)(6,0)(12,0)}\put(0,4){\bezier29(0,0)(6,0)(12,0)}}
\put(12,0){\put(1,16){$s_{41}$}\put(0,12){\put(0,0){\bezier39(0,0)(4,0)(4,-10.3)}\put(4,-12){\bezier19(0.5,-1.5)(4,-2.5)(4,1.7)}\put(8,0){\bezier39(0,-10.3)(0,0)(4,0)}}\put(0,0){\put(0,0){\bezier19(0,0)(3,0)(5.7,0)}\put(12,0){\bezier19(-2.1,0)(-1,0)(0,0)}}\put(0,4){\put(0,0){\bezier9(0,0)(1,0)(2,0)}\put(6,0){\bezier9(-0.2,0)(0,0)(0.2,0)}\put(12,0){\bezier9(-2,0)(-1,0)(0,0)}}\put(0,8){\put(0,0){\bezier9(0,0)(1,0)(1.5,0)}\put(6,0){\bezier9(-0.7,0)(0,0)(0.7,0)}\put(12,0){\bezier9(-1.5,0)(-1,0)(0,0)}}}
\put(24,0){\put(1,16){$s_{41}$}\put(0,12){\put(0,0){\bezier39(0,0)(4,0)(4,-10.3)}\put(4,-12){\bezier19(0.5,-1.5)(4,-2.5)(4,1.7)}\put(8,0){\bezier39(0,-10.3)(0,0)(4,0)}}\put(0,0){\put(0,0){\bezier19(0,0)(3,0)(5.7,0)}\put(12,0){\bezier19(-2.1,0)(-1,0)(0,0)}}\put(0,4){\put(0,0){\bezier9(0,0)(1,0)(2,0)}\put(6,0){\bezier9(-0.2,0)(0,0)(0.2,0)}\put(12,0){\bezier9(-2,0)(-1,0)(0,0)}}\put(0,8){\put(0,0){\bezier9(0,0)(1,0)(1.5,0)}\put(6,0){\bezier9(-0.7,0)(0,0)(0.7,0)}\put(12,0){\bezier9(-1.5,0)(-1,0)(0,0)}}}
\put(36,0){\put(1,16){$s_{43}^{-1}$}\put(0,12){\put(0,0){\bezier19(0,0)(4,0)(4,-2.3)}\put(4,-4){\bezier19(0,1.7)(0,-2.5)(3.5,-1.5)}\put(8,0){\bezier19(0,-2.3)(0,0)(4,0)}}\put(0,8){\put(0,0){\bezier19(0,0)(1,0)(2.1,0)}\put(12,0){\bezier19(-5.7,0)(-3,0)(0,0)}}\put(0,0){\bezier29(0,0)(6,0)(12,0)}\put(0,4){\bezier29(0,0)(6,0)(12,0)}}
\put(48,0){\put(1,16){$s_{41}^{-1}$}\put(0,12){\put(0,0){\bezier39(0,0)(4,0)(4,-10.3)}\put(4,-12){\bezier19(0,1.7)(0,-2.5)(3.5,-1.5)}\put(8,0){\bezier39(0,-10.3)(0,0)(4,0)}}\put(0,0){\put(0,0){\bezier19(0,0)(1,0)(2.1,0)}\put(12,0){\bezier19(-5.7,0)(-3,0)(0,0)}}\put(0,4){\put(0,0){\bezier9(0,0)(1,0)(2,0)}\put(6,0){\bezier9(-0.2,0)(0,0)(0.2,0)}\put(12,0){\bezier9(-2,0)(-1,0)(0,0)}}\put(0,8){\put(0,0){\bezier9(0,0)(1,0)(1.5,0)}\put(6,0){\bezier9(-0.7,0)(0,0)(0.7,0)}\put(12,0){\bezier9(-1.5,0)(-1,0)(0,0)}}}
\put(60,0){\put(1,16){$s_{43}^{-1}$}\put(0,12){\put(0,0){\bezier19(0,0)(4,0)(4,-2.3)}\put(4,-4){\bezier19(0,1.7)(0,-2.5)(3.5,-1.5)}\put(8,0){\bezier19(0,-2.3)(0,0)(4,0)}}\put(0,8){\put(0,0){\bezier19(0,0)(1,0)(2.1,0)}\put(12,0){\bezier19(-5.7,0)(-3,0)(0,0)}}\put(0,0){\bezier29(0,0)(6,0)(12,0)}\put(0,4){\bezier29(0,0)(6,0)(12,0)}}
\put(72,0){\put(1,16){$s_{31}^{-1}.$}\put(0,8){\put(0,0){\bezier29(0,0)(4,0)(4,-6.3)}\put(4,-8){\bezier19(0,1.7)(0,-2.5)(3.5,-1.5)}\put(8,0){\bezier29(0,-6.3)(0,0)(4,0)}}\put(0,0){\put(0,0){\bezier19(0,0)(1,0)(2.1,0)}\put(12,0){\bezier19(-5.7,0)(-3,0)(0,0)}}\put(0,4){\put(0,0){\bezier9(0,0)(1,0)(1.7,0)}\put(6,0){\bezier9(-0.3,0)(0,0)(0.3,0)}\put(12,0){\bezier9(-1.7,0)(-1,0)(0,0)}}\put(0,12){\bezier29(0,0)(6,0)(12,0)}}
}
\put(-43,175) {
\put(-93,16){$\gamma_4=\gamma_3\cdot
s_{32}^{-1};$}\put(-30,16){$\mathfrak{I}(\gamma_4)=$}
\put(0,0){\put(1,16){$s_{43}$}\put(0,12){\put(0,0){\bezier19(0,0)(4,0)(4,-2.3)}\put(4,-4){\bezier19(0.5,-1.5)(4,-2.5)(4,1.7)}\put(8,0){\bezier19(0,-2.3)(0,0)(4,0)}}\put(0,8){\put(0,0){\bezier19(0,0)(3,0)(5.7,0)}\put(12,0){\bezier9(-2.1,0)(-1,0)(0,0)}}\put(0,0){\bezier29(0,0)(6,0)(12,0)}\put(0,4){\bezier29(0,0)(6,0)(12,0)}}
\put(12,0){\put(1,16){$s_{41}$}\put(0,12){\put(0,0){\bezier39(0,0)(4,0)(4,-10.3)}\put(4,-12){\bezier19(0.5,-1.5)(4,-2.5)(4,1.7)}\put(8,0){\bezier39(0,-10.3)(0,0)(4,0)}}\put(0,0){\put(0,0){\bezier19(0,0)(3,0)(5.7,0)}\put(12,0){\bezier19(-2.1,0)(-1,0)(0,0)}}\put(0,4){\put(0,0){\bezier9(0,0)(1,0)(2,0)}\put(6,0){\bezier9(-0.2,0)(0,0)(0.2,0)}\put(12,0){\bezier9(-2,0)(-1,0)(0,0)}}\put(0,8){\put(0,0){\bezier9(0,0)(1,0)(1.5,0)}\put(6,0){\bezier9(-0.7,0)(0,0)(0.7,0)}\put(12,0){\bezier9(-1.5,0)(-1,0)(0,0)}}}
\put(24,0){\put(1,16){$s_{41}$}\put(0,12){\put(0,0){\bezier39(0,0)(4,0)(4,-10.3)}\put(4,-12){\bezier19(0.5,-1.5)(4,-2.5)(4,1.7)}\put(8,0){\bezier39(0,-10.3)(0,0)(4,0)}}\put(0,0){\put(0,0){\bezier19(0,0)(3,0)(5.7,0)}\put(12,0){\bezier19(-2.1,0)(-1,0)(0,0)}}\put(0,4){\put(0,0){\bezier9(0,0)(1,0)(2,0)}\put(6,0){\bezier9(-0.2,0)(0,0)(0.2,0)}\put(12,0){\bezier9(-2,0)(-1,0)(0,0)}}\put(0,8){\put(0,0){\bezier9(0,0)(1,0)(1.5,0)}\put(6,0){\bezier9(-0.7,0)(0,0)(0.7,0)}\put(12,0){\bezier9(-1.5,0)(-1,0)(0,0)}}}
\put(36,0){\put(1,16){$s_{43}^{-1}$}\put(0,12){\put(0,0){\bezier19(0,0)(4,0)(4,-2.3)}\put(4,-4){\bezier19(0,1.7)(0,-2.5)(3.5,-1.5)}\put(8,0){\bezier19(0,-2.3)(0,0)(4,0)}}\put(0,8){\put(0,0){\bezier19(0,0)(1,0)(2.1,0)}\put(12,0){\bezier19(-5.7,0)(-3,0)(0,0)}}\put(0,0){\bezier29(0,0)(6,0)(12,0)}\put(0,4){\bezier29(0,0)(6,0)(12,0)}}
\put(48,0){\put(1,16){$s_{41}^{-1}$}\put(0,12){\put(0,0){\bezier39(0,0)(4,0)(4,-10.3)}\put(4,-12){\bezier19(0,1.7)(0,-2.5)(3.5,-1.5)}\put(8,0){\bezier39(0,-10.3)(0,0)(4,0)}}\put(0,0){\put(0,0){\bezier19(0,0)(1,0)(2.1,0)}\put(12,0){\bezier19(-5.7,0)(-3,0)(0,0)}}\put(0,4){\put(0,0){\bezier9(0,0)(1,0)(2,0)}\put(6,0){\bezier9(-0.2,0)(0,0)(0.2,0)}\put(12,0){\bezier9(-2,0)(-1,0)(0,0)}}\put(0,8){\put(0,0){\bezier9(0,0)(1,0)(1.5,0)}\put(6,0){\bezier9(-0.7,0)(0,0)(0.7,0)}\put(12,0){\bezier9(-1.5,0)(-1,0)(0,0)}}}
\put(60,0){\put(1,16){$s_{43}^{-1}$}\put(0,12){\put(0,0){\bezier19(0,0)(4,0)(4,-2.3)}\put(4,-4){\bezier19(0,1.7)(0,-2.5)(3.5,-1.5)}\put(8,0){\bezier19(0,-2.3)(0,0)(4,0)}}\put(0,8){\put(0,0){\bezier19(0,0)(1,0)(2.1,0)}\put(12,0){\bezier19(-5.7,0)(-3,0)(0,0)}}\put(0,0){\bezier29(0,0)(6,0)(12,0)}\put(0,4){\bezier29(0,0)(6,0)(12,0)}}
\put(72,0){\put(1,16){$s_{31}^{-1}$}\put(0,8){\put(0,0){\bezier29(0,0)(4,0)(4,-6.3)}\put(4,-8){\bezier19(0,1.7)(0,-2.5)(3.5,-1.5)}\put(8,0){\bezier29(0,-6.3)(0,0)(4,0)}}\put(0,0){\put(0,0){\bezier19(0,0)(1,0)(2.1,0)}\put(12,0){\bezier19(-5.7,0)(-3,0)(0,0)}}\put(0,4){\put(0,0){\bezier9(0,0)(1,0)(1.7,0)}\put(6,0){\bezier9(-0.3,0)(0,0)(0.3,0)}\put(12,0){\bezier9(-1.7,0)(-1,0)(0,0)}}\put(0,12){\bezier29(0,0)(6,0)(12,0)}}
\put(84,0){\put(1,16){$s_{32}^{-1}.$}\put(0,8){\put(0,0){\bezier19(0,0)(4,0)(4,-2.3)}\put(4,-4){\bezier19(0,1.7)(0,-2.5)(3.5,-1.5)}\put(8,0){\bezier19(0,-2.3)(0,0)(4,0)}}\put(0,4){\put(0,0){\bezier19(0,0)(1,0)(2.1,0)}\put(12,0){\bezier19(-5.7,0)(-3,0)(0,0)}}\put(0,12){\bezier29(0,0)(6,0)(12,0)}\put(0,0){\bezier29(0,0)(6,0)(12,0)}}
}
\put(-43,140) {
\put(-93,16){$\gamma_5=\gamma_4\cdot
s_{42};$}\put(-30,16){$\mathfrak{I}(\gamma_5)=$}
\put(0,0){\put(1,16){$s_{43}$}\put(0,12){\put(0,0){\bezier19(0,0)(4,0)(4,-2.3)}\put(4,-4){\bezier19(0.5,-1.5)(4,-2.5)(4,1.7)}\put(8,0){\bezier19(0,-2.3)(0,0)(4,0)}}\put(0,8){\put(0,0){\bezier19(0,0)(3,0)(5.7,0)}\put(12,0){\bezier9(-2.1,0)(-1,0)(0,0)}}\put(0,0){\bezier29(0,0)(6,0)(12,0)}\put(0,4){\bezier29(0,0)(6,0)(12,0)}}
\put(12,0){\put(1,16){$s_{41}$}\put(0,12){\put(0,0){\bezier39(0,0)(4,0)(4,-10.3)}\put(4,-12){\bezier19(0.5,-1.5)(4,-2.5)(4,1.7)}\put(8,0){\bezier39(0,-10.3)(0,0)(4,0)}}\put(0,0){\put(0,0){\bezier19(0,0)(3,0)(5.7,0)}\put(12,0){\bezier19(-2.1,0)(-1,0)(0,0)}}\put(0,4){\put(0,0){\bezier9(0,0)(1,0)(2,0)}\put(6,0){\bezier9(-0.2,0)(0,0)(0.2,0)}\put(12,0){\bezier9(-2,0)(-1,0)(0,0)}}\put(0,8){\put(0,0){\bezier9(0,0)(1,0)(1.5,0)}\put(6,0){\bezier9(-0.7,0)(0,0)(0.7,0)}\put(12,0){\bezier9(-1.5,0)(-1,0)(0,0)}}}
\put(24,0){\put(1,16){$s_{41}$}\put(0,12){\put(0,0){\bezier39(0,0)(4,0)(4,-10.3)}\put(4,-12){\bezier19(0.5,-1.5)(4,-2.5)(4,1.7)}\put(8,0){\bezier39(0,-10.3)(0,0)(4,0)}}\put(0,0){\put(0,0){\bezier19(0,0)(3,0)(5.7,0)}\put(12,0){\bezier19(-2.1,0)(-1,0)(0,0)}}\put(0,4){\put(0,0){\bezier9(0,0)(1,0)(2,0)}\put(6,0){\bezier9(-0.2,0)(0,0)(0.2,0)}\put(12,0){\bezier9(-2,0)(-1,0)(0,0)}}\put(0,8){\put(0,0){\bezier9(0,0)(1,0)(1.5,0)}\put(6,0){\bezier9(-0.7,0)(0,0)(0.7,0)}\put(12,0){\bezier9(-1.5,0)(-1,0)(0,0)}}}
\put(36,0){\put(1,16){$s_{43}^{-1}$}\put(0,12){\put(0,0){\bezier19(0,0)(4,0)(4,-2.3)}\put(4,-4){\bezier19(0,1.7)(0,-2.5)(3.5,-1.5)}\put(8,0){\bezier19(0,-2.3)(0,0)(4,0)}}\put(0,8){\put(0,0){\bezier19(0,0)(1,0)(2.1,0)}\put(12,0){\bezier19(-5.7,0)(-3,0)(0,0)}}\put(0,0){\bezier29(0,0)(6,0)(12,0)}\put(0,4){\bezier29(0,0)(6,0)(12,0)}}
\put(48,0){\put(1,16){$s_{41}^{-1}$}\put(0,12){\put(0,0){\bezier39(0,0)(4,0)(4,-10.3)}\put(4,-12){\bezier19(0,1.7)(0,-2.5)(3.5,-1.5)}\put(8,0){\bezier39(0,-10.3)(0,0)(4,0)}}\put(0,0){\put(0,0){\bezier19(0,0)(1,0)(2.1,0)}\put(12,0){\bezier19(-5.7,0)(-3,0)(0,0)}}\put(0,4){\put(0,0){\bezier9(0,0)(1,0)(2,0)}\put(6,0){\bezier9(-0.2,0)(0,0)(0.2,0)}\put(12,0){\bezier9(-2,0)(-1,0)(0,0)}}\put(0,8){\put(0,0){\bezier9(0,0)(1,0)(1.5,0)}\put(6,0){\bezier9(-0.7,0)(0,0)(0.7,0)}\put(12,0){\bezier9(-1.5,0)(-1,0)(0,0)}}}
\put(60,0){\put(1,16){$s_{42}$}\put(0,12){\put(0,0){\bezier29(0,0)(4,0)(4,-6.3)}\put(4,-8){\bezier19(0.5,-1.5)(4,-2.5)(4,1.7)}\put(8,0){\bezier29(0,-6.3)(0,0)(4,0)}}\put(0,4){\put(0,0){\bezier19(0,0)(3,0)(5.7,0)}\put(12,0){\bezier19(-2.1,0)(-1,0)(0,0)}}\put(0,8){\put(0,0){\bezier9(0,0)(1,0)(1.7,0)}\put(6,0){\bezier9(-0.3,0)(0,0)(0.3,0)}\put(12,0){\bezier9(-1.7,0)(-1,0)(0,0)}}\put(0,0){\bezier29(0,0)(6,0)(12,0)}}
\put(72,0){\put(1,16){$s_{43}^{-1}$}\put(0,12){\put(0,0){\bezier19(0,0)(4,0)(4,-2.3)}\put(4,-4){\bezier19(0,1.7)(0,-2.5)(3.5,-1.5)}\put(8,0){\bezier19(0,-2.3)(0,0)(4,0)}}\put(0,8){\put(0,0){\bezier19(0,0)(1,0)(2.1,0)}\put(12,0){\bezier19(-5.7,0)(-3,0)(0,0)}}\put(0,0){\bezier29(0,0)(6,0)(12,0)}\put(0,4){\bezier29(0,0)(6,0)(12,0)}}
\put(84,0){\put(1,16){$s_{31}^{-1}$}\put(0,8){\put(0,0){\bezier29(0,0)(4,0)(4,-6.3)}\put(4,-8){\bezier19(0,1.7)(0,-2.5)(3.5,-1.5)}\put(8,0){\bezier29(0,-6.3)(0,0)(4,0)}}\put(0,0){\put(0,0){\bezier19(0,0)(1,0)(2.1,0)}\put(12,0){\bezier19(-5.7,0)(-3,0)(0,0)}}\put(0,4){\put(0,0){\bezier9(0,0)(1,0)(1.7,0)}\put(6,0){\bezier9(-0.3,0)(0,0)(0.3,0)}\put(12,0){\bezier9(-1.7,0)(-1,0)(0,0)}}\put(0,12){\bezier29(0,0)(6,0)(12,0)}}
\put(96,0){\put(1,16){$s_{32}^{-1}.$}\put(0,8){\put(0,0){\bezier19(0,0)(4,0)(4,-2.3)}\put(4,-4){\bezier19(0,1.7)(0,-2.5)(3.5,-1.5)}\put(8,0){\bezier19(0,-2.3)(0,0)(4,0)}}\put(0,4){\put(0,0){\bezier19(0,0)(1,0)(2.1,0)}\put(12,0){\bezier19(-5.7,0)(-3,0)(0,0)}}\put(0,12){\bezier29(0,0)(6,0)(12,0)}\put(0,0){\bezier29(0,0)(6,0)(12,0)}}
}
\put(-43,105) {
\put(-93,16){$\gamma_6=\gamma_5\cdot
s_{21};$}\put(-30,16){$\mathfrak{I}(\gamma_6)=$}
\put(0,0){\put(1,16){$s_{43}$}\put(0,12){\put(0,0){\bezier19(0,0)(4,0)(4,-2.3)}\put(4,-4){\bezier19(0.5,-1.5)(4,-2.5)(4,1.7)}\put(8,0){\bezier19(0,-2.3)(0,0)(4,0)}}\put(0,8){\put(0,0){\bezier19(0,0)(3,0)(5.7,0)}\put(12,0){\bezier9(-2.1,0)(-1,0)(0,0)}}\put(0,0){\bezier29(0,0)(6,0)(12,0)}\put(0,4){\bezier29(0,0)(6,0)(12,0)}}
\put(12,0){\put(1,16){$s_{41}$}\put(0,12){\put(0,0){\bezier39(0,0)(4,0)(4,-10.3)}\put(4,-12){\bezier19(0.5,-1.5)(4,-2.5)(4,1.7)}\put(8,0){\bezier39(0,-10.3)(0,0)(4,0)}}\put(0,0){\put(0,0){\bezier19(0,0)(3,0)(5.7,0)}\put(12,0){\bezier19(-2.1,0)(-1,0)(0,0)}}\put(0,4){\put(0,0){\bezier9(0,0)(1,0)(2,0)}\put(6,0){\bezier9(-0.2,0)(0,0)(0.2,0)}\put(12,0){\bezier9(-2,0)(-1,0)(0,0)}}\put(0,8){\put(0,0){\bezier9(0,0)(1,0)(1.5,0)}\put(6,0){\bezier9(-0.7,0)(0,0)(0.7,0)}\put(12,0){\bezier9(-1.5,0)(-1,0)(0,0)}}}
\put(24,0){\put(1,16){$s_{41}$}\put(0,12){\put(0,0){\bezier39(0,0)(4,0)(4,-10.3)}\put(4,-12){\bezier19(0.5,-1.5)(4,-2.5)(4,1.7)}\put(8,0){\bezier39(0,-10.3)(0,0)(4,0)}}\put(0,0){\put(0,0){\bezier19(0,0)(3,0)(5.7,0)}\put(12,0){\bezier19(-2.1,0)(-1,0)(0,0)}}\put(0,4){\put(0,0){\bezier9(0,0)(1,0)(2,0)}\put(6,0){\bezier9(-0.2,0)(0,0)(0.2,0)}\put(12,0){\bezier9(-2,0)(-1,0)(0,0)}}\put(0,8){\put(0,0){\bezier9(0,0)(1,0)(1.5,0)}\put(6,0){\bezier9(-0.7,0)(0,0)(0.7,0)}\put(12,0){\bezier9(-1.5,0)(-1,0)(0,0)}}}
\put(36,0){\put(1,16){$s_{43}^{-1}$}\put(0,12){\put(0,0){\bezier19(0,0)(4,0)(4,-2.3)}\put(4,-4){\bezier19(0,1.7)(0,-2.5)(3.5,-1.5)}\put(8,0){\bezier19(0,-2.3)(0,0)(4,0)}}\put(0,8){\put(0,0){\bezier19(0,0)(1,0)(2.1,0)}\put(12,0){\bezier19(-5.7,0)(-3,0)(0,0)}}\put(0,0){\bezier29(0,0)(6,0)(12,0)}\put(0,4){\bezier29(0,0)(6,0)(12,0)}}
\put(48,0){\put(1,16){$s_{41}^{-1}$}\put(0,12){\put(0,0){\bezier39(0,0)(4,0)(4,-10.3)}\put(4,-12){\bezier19(0,1.7)(0,-2.5)(3.5,-1.5)}\put(8,0){\bezier39(0,-10.3)(0,0)(4,0)}}\put(0,0){\put(0,0){\bezier19(0,0)(1,0)(2.1,0)}\put(12,0){\bezier19(-5.7,0)(-3,0)(0,0)}}\put(0,4){\put(0,0){\bezier9(0,0)(1,0)(2,0)}\put(6,0){\bezier9(-0.2,0)(0,0)(0.2,0)}\put(12,0){\bezier9(-2,0)(-1,0)(0,0)}}\put(0,8){\put(0,0){\bezier9(0,0)(1,0)(1.5,0)}\put(6,0){\bezier9(-0.7,0)(0,0)(0.7,0)}\put(12,0){\bezier9(-1.5,0)(-1,0)(0,0)}}}
\put(60,0){\put(1,16){$s_{42}$}\put(0,12){\put(0,0){\bezier29(0,0)(4,0)(4,-6.3)}\put(4,-8){\bezier19(0.5,-1.5)(4,-2.5)(4,1.7)}\put(8,0){\bezier29(0,-6.3)(0,0)(4,0)}}\put(0,4){\put(0,0){\bezier19(0,0)(3,0)(5.7,0)}\put(12,0){\bezier19(-2.1,0)(-1,0)(0,0)}}\put(0,8){\put(0,0){\bezier9(0,0)(1,0)(1.7,0)}\put(6,0){\bezier9(-0.3,0)(0,0)(0.3,0)}\put(12,0){\bezier9(-1.7,0)(-1,0)(0,0)}}\put(0,0){\bezier29(0,0)(6,0)(12,0)}}
\put(72,0){\put(1,16){$s_{43}^{-1}$}\put(0,12){\put(0,0){\bezier19(0,0)(4,0)(4,-2.3)}\put(4,-4){\bezier19(0,1.7)(0,-2.5)(3.5,-1.5)}\put(8,0){\bezier19(0,-2.3)(0,0)(4,0)}}\put(0,8){\put(0,0){\bezier19(0,0)(1,0)(2.1,0)}\put(12,0){\bezier19(-5.7,0)(-3,0)(0,0)}}\put(0,0){\bezier29(0,0)(6,0)(12,0)}\put(0,4){\bezier29(0,0)(6,0)(12,0)}}
\put(84,0){\put(1,16){$s_{31}^{-1}$}\put(0,8){\put(0,0){\bezier29(0,0)(4,0)(4,-6.3)}\put(4,-8){\bezier19(0,1.7)(0,-2.5)(3.5,-1.5)}\put(8,0){\bezier29(0,-6.3)(0,0)(4,0)}}\put(0,0){\put(0,0){\bezier19(0,0)(1,0)(2.1,0)}\put(12,0){\bezier19(-5.7,0)(-3,0)(0,0)}}\put(0,4){\put(0,0){\bezier9(0,0)(1,0)(1.7,0)}\put(6,0){\bezier9(-0.3,0)(0,0)(0.3,0)}\put(12,0){\bezier9(-1.7,0)(-1,0)(0,0)}}\put(0,12){\bezier29(0,0)(6,0)(12,0)}}
\put(96,0){\put(1,16){$s_{32}^{-1}$}\put(0,8){\put(0,0){\bezier19(0,0)(4,0)(4,-2.3)}\put(4,-4){\bezier19(0,1.7)(0,-2.5)(3.5,-1.5)}\put(8,0){\bezier19(0,-2.3)(0,0)(4,0)}}\put(0,4){\put(0,0){\bezier19(0,0)(1,0)(2.1,0)}\put(12,0){\bezier19(-5.7,0)(-3,0)(0,0)}}\put(0,12){\bezier29(0,0)(6,0)(12,0)}\put(0,0){\bezier29(0,0)(6,0)(12,0)}}
\put(108,0){\put(1,16){$s_{21}.$}\put(0,4){\put(0,0){\bezier19(0,0)(4,0)(4,-2.3)}\put(4,-4){\bezier19(0.5,-1.5)(4,-2.5)(4,1.7)}\put(8,0){\bezier19(0,-2.3)(0,0)(4,0)}}\put(0,0){\put(0,0){\bezier19(0,0)(3,0)(5.7,0)}\put(12,0){\bezier9(-2.1,0)(-1,0)(0,0)}}\put(0,8){\bezier29(0,0)(6,0)(12,0)}\put(0,12){\bezier29(0,0)(6,0)(12,0)}}
}
\put(-43,70) {
\put(-93,16){$\gamma_7=\gamma_6\cdot
s_{32};$}\put(-30,16){$\mathfrak{I}(\gamma_7)=$}
\put(0,0){\put(1,16){$s_{43}$}\put(0,12){\put(0,0){\bezier19(0,0)(4,0)(4,-2.3)}\put(4,-4){\bezier19(0.5,-1.5)(4,-2.5)(4,1.7)}\put(8,0){\bezier19(0,-2.3)(0,0)(4,0)}}\put(0,8){\put(0,0){\bezier19(0,0)(3,0)(5.7,0)}\put(12,0){\bezier9(-2.1,0)(-1,0)(0,0)}}\put(0,0){\bezier29(0,0)(6,0)(12,0)}\put(0,4){\bezier29(0,0)(6,0)(12,0)}}
\put(12,0){\put(1,16){$s_{41}$}\put(0,12){\put(0,0){\bezier39(0,0)(4,0)(4,-10.3)}\put(4,-12){\bezier19(0.5,-1.5)(4,-2.5)(4,1.7)}\put(8,0){\bezier39(0,-10.3)(0,0)(4,0)}}\put(0,0){\put(0,0){\bezier19(0,0)(3,0)(5.7,0)}\put(12,0){\bezier19(-2.1,0)(-1,0)(0,0)}}\put(0,4){\put(0,0){\bezier9(0,0)(1,0)(2,0)}\put(6,0){\bezier9(-0.2,0)(0,0)(0.2,0)}\put(12,0){\bezier9(-2,0)(-1,0)(0,0)}}\put(0,8){\put(0,0){\bezier9(0,0)(1,0)(1.5,0)}\put(6,0){\bezier9(-0.7,0)(0,0)(0.7,0)}\put(12,0){\bezier9(-1.5,0)(-1,0)(0,0)}}}
\put(24,0){\put(1,16){$s_{41}$}\put(0,12){\put(0,0){\bezier39(0,0)(4,0)(4,-10.3)}\put(4,-12){\bezier19(0.5,-1.5)(4,-2.5)(4,1.7)}\put(8,0){\bezier39(0,-10.3)(0,0)(4,0)}}\put(0,0){\put(0,0){\bezier19(0,0)(3,0)(5.7,0)}\put(12,0){\bezier19(-2.1,0)(-1,0)(0,0)}}\put(0,4){\put(0,0){\bezier9(0,0)(1,0)(2,0)}\put(6,0){\bezier9(-0.2,0)(0,0)(0.2,0)}\put(12,0){\bezier9(-2,0)(-1,0)(0,0)}}\put(0,8){\put(0,0){\bezier9(0,0)(1,0)(1.5,0)}\put(6,0){\bezier9(-0.7,0)(0,0)(0.7,0)}\put(12,0){\bezier9(-1.5,0)(-1,0)(0,0)}}}
\put(36,0){\put(1,16){$s_{43}^{-1}$}\put(0,12){\put(0,0){\bezier19(0,0)(4,0)(4,-2.3)}\put(4,-4){\bezier19(0,1.7)(0,-2.5)(3.5,-1.5)}\put(8,0){\bezier19(0,-2.3)(0,0)(4,0)}}\put(0,8){\put(0,0){\bezier19(0,0)(1,0)(2.1,0)}\put(12,0){\bezier19(-5.7,0)(-3,0)(0,0)}}\put(0,0){\bezier29(0,0)(6,0)(12,0)}\put(0,4){\bezier29(0,0)(6,0)(12,0)}}
\put(48,0){\put(1,16){$s_{41}^{-1}$}\put(0,12){\put(0,0){\bezier39(0,0)(4,0)(4,-10.3)}\put(4,-12){\bezier19(0,1.7)(0,-2.5)(3.5,-1.5)}\put(8,0){\bezier39(0,-10.3)(0,0)(4,0)}}\put(0,0){\put(0,0){\bezier19(0,0)(1,0)(2.1,0)}\put(12,0){\bezier19(-5.7,0)(-3,0)(0,0)}}\put(0,4){\put(0,0){\bezier9(0,0)(1,0)(2,0)}\put(6,0){\bezier9(-0.2,0)(0,0)(0.2,0)}\put(12,0){\bezier9(-2,0)(-1,0)(0,0)}}\put(0,8){\put(0,0){\bezier9(0,0)(1,0)(1.5,0)}\put(6,0){\bezier9(-0.7,0)(0,0)(0.7,0)}\put(12,0){\bezier9(-1.5,0)(-1,0)(0,0)}}}
\put(60,0){\put(1,16){$s_{42}$}\put(0,12){\put(0,0){\bezier29(0,0)(4,0)(4,-6.3)}\put(4,-8){\bezier19(0.5,-1.5)(4,-2.5)(4,1.7)}\put(8,0){\bezier29(0,-6.3)(0,0)(4,0)}}\put(0,4){\put(0,0){\bezier19(0,0)(3,0)(5.7,0)}\put(12,0){\bezier19(-2.1,0)(-1,0)(0,0)}}\put(0,8){\put(0,0){\bezier9(0,0)(1,0)(1.7,0)}\put(6,0){\bezier9(-0.3,0)(0,0)(0.3,0)}\put(12,0){\bezier9(-1.7,0)(-1,0)(0,0)}}\put(0,0){\bezier29(0,0)(6,0)(12,0)}}
\put(72,0){\put(1,16){$s_{43}^{-1}$}\put(0,12){\put(0,0){\bezier19(0,0)(4,0)(4,-2.3)}\put(4,-4){\bezier19(0,1.7)(0,-2.5)(3.5,-1.5)}\put(8,0){\bezier19(0,-2.3)(0,0)(4,0)}}\put(0,8){\put(0,0){\bezier19(0,0)(1,0)(2.1,0)}\put(12,0){\bezier19(-5.7,0)(-3,0)(0,0)}}\put(0,0){\bezier29(0,0)(6,0)(12,0)}\put(0,4){\bezier29(0,0)(6,0)(12,0)}}
\put(84,0){\put(1,16){$s_{31}^{-1}$}\put(0,8){\put(0,0){\bezier29(0,0)(4,0)(4,-6.3)}\put(4,-8){\bezier19(0,1.7)(0,-2.5)(3.5,-1.5)}\put(8,0){\bezier29(0,-6.3)(0,0)(4,0)}}\put(0,0){\put(0,0){\bezier19(0,0)(1,0)(2.1,0)}\put(12,0){\bezier19(-5.7,0)(-3,0)(0,0)}}\put(0,4){\put(0,0){\bezier9(0,0)(1,0)(1.7,0)}\put(6,0){\bezier9(-0.3,0)(0,0)(0.3,0)}\put(12,0){\bezier9(-1.7,0)(-1,0)(0,0)}}\put(0,12){\bezier29(0,0)(6,0)(12,0)}}
\put(96,0){\put(1,16){$s_{32}^{-1}$}\put(0,8){\put(0,0){\bezier19(0,0)(4,0)(4,-2.3)}\put(4,-4){\bezier19(0,1.7)(0,-2.5)(3.5,-1.5)}\put(8,0){\bezier19(0,-2.3)(0,0)(4,0)}}\put(0,4){\put(0,0){\bezier19(0,0)(1,0)(2.1,0)}\put(12,0){\bezier19(-5.7,0)(-3,0)(0,0)}}\put(0,12){\bezier29(0,0)(6,0)(12,0)}\put(0,0){\bezier29(0,0)(6,0)(12,0)}}
\put(108,0){\put(1,16){$s_{31}^{-1}$}\put(0,8){\put(0,0){\bezier29(0,0)(4,0)(4,-6.3)}\put(4,-8){\bezier19(0,1.7)(0,-2.5)(3.5,-1.5)}\put(8,0){\bezier29(0,-6.3)(0,0)(4,0)}}\put(0,0){\put(0,0){\bezier19(0,0)(1,0)(2.1,0)}\put(12,0){\bezier19(-5.7,0)(-3,0)(0,0)}}\put(0,4){\put(0,0){\bezier9(0,0)(1,0)(1.7,0)}\put(6,0){\bezier9(-0.3,0)(0,0)(0.3,0)}\put(12,0){\bezier9(-1.7,0)(-1,0)(0,0)}}\put(0,12){\bezier29(0,0)(6,0)(12,0)}}
\put(120,0){\put(1,16){$s_{32}$}\put(0,8){\put(0,0){\bezier19(0,0)(4,0)(4,-2.3)}\put(4,-4){\bezier19(0.5,-1.5)(4,-2.5)(4,1.7)}\put(8,0){\bezier19(0,-2.3)(0,0)(4,0)}}\put(0,4){\put(0,0){\bezier19(0,0)(3,0)(5.7,0)}\put(12,0){\bezier9(-2.1,0)(-1,0)(0,0)}}\put(0,12){\bezier29(0,0)(6,0)(12,0)}\put(0,0){\bezier29(0,0)(6,0)(12,0)}}
\put(132,0){\put(1,16){$s_{31}$}\put(0,8){\put(0,0){\bezier29(0,0)(4,0)(4,-6.3)}\put(4,-8){\bezier19(0.5,-1.5)(4,-2.5)(4,1.7)}\put(8,0){\bezier29(0,-6.3)(0,0)(4,0)}}\put(0,0){\put(0,0){\bezier19(0,0)(3,0)(5.7,0)}\put(12,0){\bezier19(-2.1,0)(-1,0)(0,0)}}\put(0,4){\put(0,0){\bezier9(0,0)(1,0)(1.7,0)}\put(6,0){\bezier9(-0.3,0)(0,0)(0.3,0)}\put(12,0){\bezier9(-1.7,0)(-1,0)(0,0)}}\put(0,12){\bezier29(0,0)(6,0)(12,0)}}
\put(144,0){\put(1,16){$s_{21}.$}\put(0,4){\put(0,0){\bezier19(0,0)(4,0)(4,-2.3)}\put(4,-4){\bezier19(0.5,-1.5)(4,-2.5)(4,1.7)}\put(8,0){\bezier19(0,-2.3)(0,0)(4,0)}}\put(0,0){\put(0,0){\bezier19(0,0)(3,0)(5.7,0)}\put(12,0){\bezier9(-2.1,0)(-1,0)(0,0)}}\put(0,8){\bezier29(0,0)(6,0)(12,0)}\put(0,12){\bezier29(0,0)(6,0)(12,0)}}
}
\put(-43,35) {
\put(-93,16){$\gamma_8=\gamma_7\cdot
s_{41}^{-1};$}\put(-30,16){$\mathfrak{I}(\gamma_8)=$}
\put(0,0){\put(1,16){$s_{43}$}\put(0,12){\put(0,0){\bezier19(0,0)(4,0)(4,-2.3)}\put(4,-4){\bezier19(0.5,-1.5)(4,-2.5)(4,1.7)}\put(8,0){\bezier19(0,-2.3)(0,0)(4,0)}}\put(0,8){\put(0,0){\bezier19(0,0)(3,0)(5.7,0)}\put(12,0){\bezier9(-2.1,0)(-1,0)(0,0)}}\put(0,0){\bezier29(0,0)(6,0)(12,0)}\put(0,4){\bezier29(0,0)(6,0)(12,0)}}
\put(12,0){\put(1,16){$s_{41}$}\put(0,12){\put(0,0){\bezier39(0,0)(4,0)(4,-10.3)}\put(4,-12){\bezier19(0.5,-1.5)(4,-2.5)(4,1.7)}\put(8,0){\bezier39(0,-10.3)(0,0)(4,0)}}\put(0,0){\put(0,0){\bezier19(0,0)(3,0)(5.7,0)}\put(12,0){\bezier19(-2.1,0)(-1,0)(0,0)}}\put(0,4){\put(0,0){\bezier9(0,0)(1,0)(2,0)}\put(6,0){\bezier9(-0.2,0)(0,0)(0.2,0)}\put(12,0){\bezier9(-2,0)(-1,0)(0,0)}}\put(0,8){\put(0,0){\bezier9(0,0)(1,0)(1.5,0)}\put(6,0){\bezier9(-0.7,0)(0,0)(0.7,0)}\put(12,0){\bezier9(-1.5,0)(-1,0)(0,0)}}}
\put(24,0){\put(1,16){$s_{41}$}\put(0,12){\put(0,0){\bezier39(0,0)(4,0)(4,-10.3)}\put(4,-12){\bezier19(0.5,-1.5)(4,-2.5)(4,1.7)}\put(8,0){\bezier39(0,-10.3)(0,0)(4,0)}}\put(0,0){\put(0,0){\bezier19(0,0)(3,0)(5.7,0)}\put(12,0){\bezier19(-2.1,0)(-1,0)(0,0)}}\put(0,4){\put(0,0){\bezier9(0,0)(1,0)(2,0)}\put(6,0){\bezier9(-0.2,0)(0,0)(0.2,0)}\put(12,0){\bezier9(-2,0)(-1,0)(0,0)}}\put(0,8){\put(0,0){\bezier9(0,0)(1,0)(1.5,0)}\put(6,0){\bezier9(-0.7,0)(0,0)(0.7,0)}\put(12,0){\bezier9(-1.5,0)(-1,0)(0,0)}}}
\put(36,0){\put(1,16){$s_{43}^{-1}$}\put(0,12){\put(0,0){\bezier19(0,0)(4,0)(4,-2.3)}\put(4,-4){\bezier19(0,1.7)(0,-2.5)(3.5,-1.5)}\put(8,0){\bezier19(0,-2.3)(0,0)(4,0)}}\put(0,8){\put(0,0){\bezier19(0,0)(1,0)(2.1,0)}\put(12,0){\bezier19(-5.7,0)(-3,0)(0,0)}}\put(0,0){\bezier29(0,0)(6,0)(12,0)}\put(0,4){\bezier29(0,0)(6,0)(12,0)}}
\put(48,0){\put(1,16){$s_{41}^{-1}$}\put(0,12){\put(0,0){\bezier39(0,0)(4,0)(4,-10.3)}\put(4,-12){\bezier19(0,1.7)(0,-2.5)(3.5,-1.5)}\put(8,0){\bezier39(0,-10.3)(0,0)(4,0)}}\put(0,0){\put(0,0){\bezier19(0,0)(1,0)(2.1,0)}\put(12,0){\bezier19(-5.7,0)(-3,0)(0,0)}}\put(0,4){\put(0,0){\bezier9(0,0)(1,0)(2,0)}\put(6,0){\bezier9(-0.2,0)(0,0)(0.2,0)}\put(12,0){\bezier9(-2,0)(-1,0)(0,0)}}\put(0,8){\put(0,0){\bezier9(0,0)(1,0)(1.5,0)}\put(6,0){\bezier9(-0.7,0)(0,0)(0.7,0)}\put(12,0){\bezier9(-1.5,0)(-1,0)(0,0)}}}
\put(60,0){\put(1,16){$s_{42}$}\put(0,12){\put(0,0){\bezier29(0,0)(4,0)(4,-6.3)}\put(4,-8){\bezier19(0.5,-1.5)(4,-2.5)(4,1.7)}\put(8,0){\bezier29(0,-6.3)(0,0)(4,0)}}\put(0,4){\put(0,0){\bezier19(0,0)(3,0)(5.7,0)}\put(12,0){\bezier19(-2.1,0)(-1,0)(0,0)}}\put(0,8){\put(0,0){\bezier9(0,0)(1,0)(1.7,0)}\put(6,0){\bezier9(-0.3,0)(0,0)(0.3,0)}\put(12,0){\bezier9(-1.7,0)(-1,0)(0,0)}}\put(0,0){\bezier29(0,0)(6,0)(12,0)}}
\put(72,0){\put(1,16){$s_{41}^{-1}$}\put(0,12){\put(0,0){\bezier39(0,0)(4,0)(4,-10.3)}\put(4,-12){\bezier19(0,1.7)(0,-2.5)(3.5,-1.5)}\put(8,0){\bezier39(0,-10.3)(0,0)(4,0)}}\put(0,0){\put(0,0){\bezier19(0,0)(1,0)(2.1,0)}\put(12,0){\bezier19(-5.7,0)(-3,0)(0,0)}}\put(0,4){\put(0,0){\bezier9(0,0)(1,0)(2,0)}\put(6,0){\bezier9(-0.2,0)(0,0)(0.2,0)}\put(12,0){\bezier9(-2,0)(-1,0)(0,0)}}\put(0,8){\put(0,0){\bezier9(0,0)(1,0)(1.5,0)}\put(6,0){\bezier9(-0.7,0)(0,0)(0.7,0)}\put(12,0){\bezier9(-1.5,0)(-1,0)(0,0)}}}
\put(84,0){\put(1,16){$s_{42}^{-1}$}\put(0,12){\put(0,0){\bezier29(0,0)(4,0)(4,-6.3)}\put(4,-8){\bezier19(0,1.7)(0,-2.5)(3.5,-1.5)}\put(8,0){\bezier29(0,-6.3)(0,0)(4,0)}}\put(0,4){\put(0,0){\bezier19(0,0)(1,0)(2.1,0)}\put(12,0){\bezier19(-5.7,0)(-3,0)(0,0)}}\put(0,8){\put(0,0){\bezier9(0,0)(1,0)(1.7,0)}\put(6,0){\bezier9(-0.3,0)(0,0)(0.3,0)}\put(12,0){\bezier9(-1.7,0)(-1,0)(0,0)}}\put(0,0){\bezier29(0,0)(6,0)(12,0)}}
\put(96,0){\put(1,16){$s_{41}^{-1}$}\put(0,12){\put(0,0){\bezier39(0,0)(4,0)(4,-10.3)}\put(4,-12){\bezier19(0,1.7)(0,-2.5)(3.5,-1.5)}\put(8,0){\bezier39(0,-10.3)(0,0)(4,0)}}\put(0,0){\put(0,0){\bezier19(0,0)(1,0)(2.1,0)}\put(12,0){\bezier19(-5.7,0)(-3,0)(0,0)}}\put(0,4){\put(0,0){\bezier9(0,0)(1,0)(2,0)}\put(6,0){\bezier9(-0.2,0)(0,0)(0.2,0)}\put(12,0){\bezier9(-2,0)(-1,0)(0,0)}}\put(0,8){\put(0,0){\bezier9(0,0)(1,0)(1.5,0)}\put(6,0){\bezier9(-0.7,0)(0,0)(0.7,0)}\put(12,0){\bezier9(-1.5,0)(-1,0)(0,0)}}}
\put(108,0){\put(1,16){$s_{42}$}\put(0,12){\put(0,0){\bezier29(0,0)(4,0)(4,-6.3)}\put(4,-8){\bezier19(0.5,-1.5)(4,-2.5)(4,1.7)}\put(8,0){\bezier29(0,-6.3)(0,0)(4,0)}}\put(0,4){\put(0,0){\bezier19(0,0)(3,0)(5.7,0)}\put(12,0){\bezier19(-2.1,0)(-1,0)(0,0)}}\put(0,8){\put(0,0){\bezier9(0,0)(1,0)(1.7,0)}\put(6,0){\bezier9(-0.3,0)(0,0)(0.3,0)}\put(12,0){\bezier9(-1.7,0)(-1,0)(0,0)}}\put(0,0){\bezier29(0,0)(6,0)(12,0)}}
\put(120,0){\put(1,16){$s_{41}$}\put(0,12){\put(0,0){\bezier39(0,0)(4,0)(4,-10.3)}\put(4,-12){\bezier19(0.5,-1.5)(4,-2.5)(4,1.7)}\put(8,0){\bezier39(0,-10.3)(0,0)(4,0)}}\put(0,0){\put(0,0){\bezier19(0,0)(3,0)(5.7,0)}\put(12,0){\bezier19(-2.1,0)(-1,0)(0,0)}}\put(0,4){\put(0,0){\bezier9(0,0)(1,0)(2,0)}\put(6,0){\bezier9(-0.2,0)(0,0)(0.2,0)}\put(12,0){\bezier9(-2,0)(-1,0)(0,0)}}\put(0,8){\put(0,0){\bezier9(0,0)(1,0)(1.5,0)}\put(6,0){\bezier9(-0.7,0)(0,0)(0.7,0)}\put(12,0){\bezier9(-1.5,0)(-1,0)(0,0)}}}
\put(132,0){\put(1,16){$s_{43}^{-1}$}\put(0,12){\put(0,0){\bezier19(0,0)(4,0)(4,-2.3)}\put(4,-4){\bezier19(0,1.7)(0,-2.5)(3.5,-1.5)}\put(8,0){\bezier19(0,-2.3)(0,0)(4,0)}}\put(0,8){\put(0,0){\bezier19(0,0)(1,0)(2.1,0)}\put(12,0){\bezier19(-5.7,0)(-3,0)(0,0)}}\put(0,0){\bezier29(0,0)(6,0)(12,0)}\put(0,4){\bezier29(0,0)(6,0)(12,0)}}
\put(144,0){\put(1,16){$s_{31}^{-1}$}\put(0,8){\put(0,0){\bezier29(0,0)(4,0)(4,-6.3)}\put(4,-8){\bezier19(0,1.7)(0,-2.5)(3.5,-1.5)}\put(8,0){\bezier29(0,-6.3)(0,0)(4,0)}}\put(0,0){\put(0,0){\bezier19(0,0)(1,0)(2.1,0)}\put(12,0){\bezier19(-5.7,0)(-3,0)(0,0)}}\put(0,4){\put(0,0){\bezier9(0,0)(1,0)(1.7,0)}\put(6,0){\bezier9(-0.3,0)(0,0)(0.3,0)}\put(12,0){\bezier9(-1.7,0)(-1,0)(0,0)}}\put(0,12){\bezier29(0,0)(6,0)(12,0)}}
\put(156,0){\put(1,16){$s_{32}^{-1}$}\put(0,8){\put(0,0){\bezier19(0,0)(4,0)(4,-2.3)}\put(4,-4){\bezier19(0,1.7)(0,-2.5)(3.5,-1.5)}\put(8,0){\bezier19(0,-2.3)(0,0)(4,0)}}\put(0,4){\put(0,0){\bezier19(0,0)(1,0)(2.1,0)}\put(12,0){\bezier19(-5.7,0)(-3,0)(0,0)}}\put(0,12){\bezier29(0,0)(6,0)(12,0)}\put(0,0){\bezier29(0,0)(6,0)(12,0)}}
\put(168,0){\put(1,16){$\cdots$}\put(0,0){\multiput(0,0)(1,4){4}{\multiput(0,0)(3,0){3}{\bezier3(0,0)(0.5,0)(0.6,0)}}}}
}
\put(-43,0) {
\put(-93,16){$\gamma_9=\gamma_8\cdot
s_{42}^{-1};$}\put(-30,16){$\mathfrak{I}(\gamma_9)=$}
\put(0,0){\put(1,16){$s_{43}$}\put(0,12){\put(0,0){\bezier19(0,0)(4,0)(4,-2.3)}\put(4,-4){\bezier19(0.5,-1.5)(4,-2.5)(4,1.7)}\put(8,0){\bezier19(0,-2.3)(0,0)(4,0)}}\put(0,8){\put(0,0){\bezier19(0,0)(3,0)(5.7,0)}\put(12,0){\bezier9(-2.1,0)(-1,0)(0,0)}}\put(0,0){\bezier29(0,0)(6,0)(12,0)}\put(0,4){\bezier29(0,0)(6,0)(12,0)}}
\put(12,0){\put(1,16){$s_{41}$}\put(0,12){\put(0,0){\bezier39(0,0)(4,0)(4,-10.3)}\put(4,-12){\bezier19(0.5,-1.5)(4,-2.5)(4,1.7)}\put(8,0){\bezier39(0,-10.3)(0,0)(4,0)}}\put(0,0){\put(0,0){\bezier19(0,0)(3,0)(5.7,0)}\put(12,0){\bezier19(-2.1,0)(-1,0)(0,0)}}\put(0,4){\put(0,0){\bezier9(0,0)(1,0)(2,0)}\put(6,0){\bezier9(-0.2,0)(0,0)(0.2,0)}\put(12,0){\bezier9(-2,0)(-1,0)(0,0)}}\put(0,8){\put(0,0){\bezier9(0,0)(1,0)(1.5,0)}\put(6,0){\bezier9(-0.7,0)(0,0)(0.7,0)}\put(12,0){\bezier9(-1.5,0)(-1,0)(0,0)}}}
\put(24,0){\put(1,16){$s_{41}$}\put(0,12){\put(0,0){\bezier39(0,0)(4,0)(4,-10.3)}\put(4,-12){\bezier19(0.5,-1.5)(4,-2.5)(4,1.7)}\put(8,0){\bezier39(0,-10.3)(0,0)(4,0)}}\put(0,0){\put(0,0){\bezier19(0,0)(3,0)(5.7,0)}\put(12,0){\bezier19(-2.1,0)(-1,0)(0,0)}}\put(0,4){\put(0,0){\bezier9(0,0)(1,0)(2,0)}\put(6,0){\bezier9(-0.2,0)(0,0)(0.2,0)}\put(12,0){\bezier9(-2,0)(-1,0)(0,0)}}\put(0,8){\put(0,0){\bezier9(0,0)(1,0)(1.5,0)}\put(6,0){\bezier9(-0.7,0)(0,0)(0.7,0)}\put(12,0){\bezier9(-1.5,0)(-1,0)(0,0)}}}
\put(36,0){\put(1,16){$s_{43}^{-1}$}\put(0,12){\put(0,0){\bezier19(0,0)(4,0)(4,-2.3)}\put(4,-4){\bezier19(0,1.7)(0,-2.5)(3.5,-1.5)}\put(8,0){\bezier19(0,-2.3)(0,0)(4,0)}}\put(0,8){\put(0,0){\bezier19(0,0)(1,0)(2.1,0)}\put(12,0){\bezier19(-5.7,0)(-3,0)(0,0)}}\put(0,0){\bezier29(0,0)(6,0)(12,0)}\put(0,4){\bezier29(0,0)(6,0)(12,0)}}
\put(48,0){\put(1,16){$s_{41}^{-1}$}\put(0,12){\put(0,0){\bezier39(0,0)(4,0)(4,-10.3)}\put(4,-12){\bezier19(0,1.7)(0,-2.5)(3.5,-1.5)}\put(8,0){\bezier39(0,-10.3)(0,0)(4,0)}}\put(0,0){\put(0,0){\bezier19(0,0)(1,0)(2.1,0)}\put(12,0){\bezier19(-5.7,0)(-3,0)(0,0)}}\put(0,4){\put(0,0){\bezier9(0,0)(1,0)(2,0)}\put(6,0){\bezier9(-0.2,0)(0,0)(0.2,0)}\put(12,0){\bezier9(-2,0)(-1,0)(0,0)}}\put(0,8){\put(0,0){\bezier9(0,0)(1,0)(1.5,0)}\put(6,0){\bezier9(-0.7,0)(0,0)(0.7,0)}\put(12,0){\bezier9(-1.5,0)(-1,0)(0,0)}}}
\put(60,0){\put(1,16){$s_{42}$}\put(0,12){\put(0,0){\bezier29(0,0)(4,0)(4,-6.3)}\put(4,-8){\bezier19(0.5,-1.5)(4,-2.5)(4,1.7)}\put(8,0){\bezier29(0,-6.3)(0,0)(4,0)}}\put(0,4){\put(0,0){\bezier19(0,0)(3,0)(5.7,0)}\put(12,0){\bezier19(-2.1,0)(-1,0)(0,0)}}\put(0,8){\put(0,0){\bezier9(0,0)(1,0)(1.7,0)}\put(6,0){\bezier9(-0.3,0)(0,0)(0.3,0)}\put(12,0){\bezier9(-1.7,0)(-1,0)(0,0)}}\put(0,0){\bezier29(0,0)(6,0)(12,0)}}
\put(72,0){\put(1,16){$s_{43}^{-1}$}\put(0,12){\put(0,0){\bezier19(0,0)(4,0)(4,-2.3)}\put(4,-4){\bezier19(0,1.7)(0,-2.5)(3.5,-1.5)}\put(8,0){\bezier19(0,-2.3)(0,0)(4,0)}}\put(0,8){\put(0,0){\bezier19(0,0)(1,0)(2.1,0)}\put(12,0){\bezier19(-5.7,0)(-3,0)(0,0)}}\put(0,0){\bezier29(0,0)(6,0)(12,0)}\put(0,4){\bezier29(0,0)(6,0)(12,0)}}
\put(84,0){\put(1,16){$s_{41}^{-1}$}\put(0,12){\put(0,0){\bezier39(0,0)(4,0)(4,-10.3)}\put(4,-12){\bezier19(0,1.7)(0,-2.5)(3.5,-1.5)}\put(8,0){\bezier39(0,-10.3)(0,0)(4,0)}}\put(0,0){\put(0,0){\bezier19(0,0)(1,0)(2.1,0)}\put(12,0){\bezier19(-5.7,0)(-3,0)(0,0)}}\put(0,4){\put(0,0){\bezier9(0,0)(1,0)(2,0)}\put(6,0){\bezier9(-0.2,0)(0,0)(0.2,0)}\put(12,0){\bezier9(-2,0)(-1,0)(0,0)}}\put(0,8){\put(0,0){\bezier9(0,0)(1,0)(1.5,0)}\put(6,0){\bezier9(-0.7,0)(0,0)(0.7,0)}\put(12,0){\bezier9(-1.5,0)(-1,0)(0,0)}}}
\put(96,0){\put(1,16){$s_{42}^{-1}$}\put(0,12){\put(0,0){\bezier29(0,0)(4,0)(4,-6.3)}\put(4,-8){\bezier19(0,1.7)(0,-2.5)(3.5,-1.5)}\put(8,0){\bezier29(0,-6.3)(0,0)(4,0)}}\put(0,4){\put(0,0){\bezier19(0,0)(1,0)(2.1,0)}\put(12,0){\bezier19(-5.7,0)(-3,0)(0,0)}}\put(0,8){\put(0,0){\bezier9(0,0)(1,0)(1.7,0)}\put(6,0){\bezier9(-0.3,0)(0,0)(0.3,0)}\put(12,0){\bezier9(-1.7,0)(-1,0)(0,0)}}\put(0,0){\bezier29(0,0)(6,0)(12,0)}}
\put(108,0){\put(1,16){$s_{41}^{-1}$}\put(0,12){\put(0,0){\bezier39(0,0)(4,0)(4,-10.3)}\put(4,-12){\bezier19(0,1.7)(0,-2.5)(3.5,-1.5)}\put(8,0){\bezier39(0,-10.3)(0,0)(4,0)}}\put(0,0){\put(0,0){\bezier19(0,0)(1,0)(2.1,0)}\put(12,0){\bezier19(-5.7,0)(-3,0)(0,0)}}\put(0,4){\put(0,0){\bezier9(0,0)(1,0)(2,0)}\put(6,0){\bezier9(-0.2,0)(0,0)(0.2,0)}\put(12,0){\bezier9(-2,0)(-1,0)(0,0)}}\put(0,8){\put(0,0){\bezier9(0,0)(1,0)(1.5,0)}\put(6,0){\bezier9(-0.7,0)(0,0)(0.7,0)}\put(12,0){\bezier9(-1.5,0)(-1,0)(0,0)}}}
\put(120,0){\put(1,16){$s_{42}^{-1}$}\put(0,12){\put(0,0){\bezier29(0,0)(4,0)(4,-6.3)}\put(4,-8){\bezier19(0,1.7)(0,-2.5)(3.5,-1.5)}\put(8,0){\bezier29(0,-6.3)(0,0)(4,0)}}\put(0,4){\put(0,0){\bezier19(0,0)(1,0)(2.1,0)}\put(12,0){\bezier19(-5.7,0)(-3,0)(0,0)}}\put(0,8){\put(0,0){\bezier9(0,0)(1,0)(1.7,0)}\put(6,0){\bezier9(-0.3,0)(0,0)(0.3,0)}\put(12,0){\bezier9(-1.7,0)(-1,0)(0,0)}}\put(0,0){\bezier29(0,0)(6,0)(12,0)}}
\put(132,0){\put(1,16){$s_{41}$}\put(0,12){\put(0,0){\bezier39(0,0)(4,0)(4,-10.3)}\put(4,-12){\bezier19(0.5,-1.5)(4,-2.5)(4,1.7)}\put(8,0){\bezier39(0,-10.3)(0,0)(4,0)}}\put(0,0){\put(0,0){\bezier19(0,0)(3,0)(5.7,0)}\put(12,0){\bezier19(-2.1,0)(-1,0)(0,0)}}\put(0,4){\put(0,0){\bezier9(0,0)(1,0)(2,0)}\put(6,0){\bezier9(-0.2,0)(0,0)(0.2,0)}\put(12,0){\bezier9(-2,0)(-1,0)(0,0)}}\put(0,8){\put(0,0){\bezier9(0,0)(1,0)(1.5,0)}\put(6,0){\bezier9(-0.7,0)(0,0)(0.7,0)}\put(12,0){\bezier9(-1.5,0)(-1,0)(0,0)}}}
\put(144,0){\put(1,16){$s_{42}$}\put(0,12){\put(0,0){\bezier29(0,0)(4,0)(4,-6.3)}\put(4,-8){\bezier19(0.5,-1.5)(4,-2.5)(4,1.7)}\put(8,0){\bezier29(0,-6.3)(0,0)(4,0)}}\put(0,4){\put(0,0){\bezier19(0,0)(3,0)(5.7,0)}\put(12,0){\bezier19(-2.1,0)(-1,0)(0,0)}}\put(0,8){\put(0,0){\bezier9(0,0)(1,0)(1.7,0)}\put(6,0){\bezier9(-0.3,0)(0,0)(0.3,0)}\put(12,0){\bezier9(-1.7,0)(-1,0)(0,0)}}\put(0,0){\bezier29(0,0)(6,0)(12,0)}}
\put(156,0){\put(1,16){$s_{41}$}\put(0,12){\put(0,0){\bezier39(0,0)(4,0)(4,-10.3)}\put(4,-12){\bezier19(0.5,-1.5)(4,-2.5)(4,1.7)}\put(8,0){\bezier39(0,-10.3)(0,0)(4,0)}}\put(0,0){\put(0,0){\bezier19(0,0)(3,0)(5.7,0)}\put(12,0){\bezier19(-2.1,0)(-1,0)(0,0)}}\put(0,4){\put(0,0){\bezier9(0,0)(1,0)(2,0)}\put(6,0){\bezier9(-0.2,0)(0,0)(0.2,0)}\put(12,0){\bezier9(-2,0)(-1,0)(0,0)}}\put(0,8){\put(0,0){\bezier9(0,0)(1,0)(1.5,0)}\put(6,0){\bezier9(-0.7,0)(0,0)(0.7,0)}\put(12,0){\bezier9(-1.5,0)(-1,0)(0,0)}}}
\put(168,0){\put(1,16){$\cdots$}\put(0,0){\multiput(0,0)(1,4){4}{\multiput(0,0)(3,0){3}{\bezier3(0,0)(0.5,0)(0.6,0)}}}}
}
} 
\end{picture}}
\caption{Geometric presentations for the Markov--Ivanovsky normal forms
(first elements in the path).}
\label{geom-path}
\end{figure}

\subsection{Corollaries about stable functionals on the braid group}
The braid group has numerous algebraic and geometric interpretations. The
Markov--Ivanovsky normal form is nicely related to some of these
interpretations, so that our results provide corresponding corollaries, which
enable us to observe numerous natural stable functionals for the random walks
on the braid group. (To describe these stable functionals briefly, we restrict
ourselves to the case of pure braid groups.)

\paragraph{The Artin representation.}
Our Theorem~3 implies that the stabilization phenomenon for the random walk on
the braid group can be observed via {\it Artin's representation\/} of $B_n$
into the automorphisms of the rank $n$ free group $F_n$ with generators
$(x_1,\ldots,x_n)$.\footnote{We remark that in the semidirect product
$P_n=F_{n-1}\rtimes P_{n-1}$ the action of the subgroup $P_{n-1}$ on the
subgroup $F_{n-1}$ yields exactly the Artin representation ``of rank $n-1$''.}
This may be described as follows. Let $x_i$ be a generator in
$(x_1,\ldots,x_n)$, and let $\gamma(x_i)$ denote the image of $x_i$ under the
automorphism that corresponds to a pure braid $\gamma$. It is well known that
the reduced word for the element $\gamma(x_i)$ has the form
$A_i(\gamma)x_iA_i^{-1}(\gamma)$ (here, $A_i(\gamma)$ is a reduced word over
the generators $x_1,\ldots,x_n$ and their inverses). It turns out that initial
(left) symbols of the Markov--Ivanovsky normal form $\mathfrak{I}(\gamma)$ for
$\gamma$ determine completely several terminal (right) letters of the word
$A_n(\gamma)$. (Equivalently, they determine initial (left) letters of the
word $A_n^{-1}(\gamma)$.) In particular, the stability of the
Markov--Ivanovsky normal form implies that for a.e. path
$\gamma_1,\gamma_2,\ldots$ of the random walk the corresponding sequence
$A_n^{-1}(\gamma_k)$ converges. Furthermore, by the symmetry arguments it
follows that for a.e. path and for each $i$ the sequence $A_i(\gamma_k)$
converges.

\begin{proclaim}{Corollary}
Let $P_n\to \operatorname{Aut}(F_n)$ be the Artin representation in its
standard form. {\rm(}In particular, for every $\gamma\in P_n$ and every $i\in
\{1,\ldots,n\}$ we have $\gamma(x_i)=A_i(\gamma)x_iA_i^{-1}(\gamma)$.\/{\rm)}
Then for a.e. path $\{\gamma_k\}_{k\in\Z_+}$ of the right random walk on $P_n$
{\rm(}with respect to a nondegenerate distribution $\mu$\/{\rm)} the sequence of
words $A_i^{-1}(\gamma_k)$ converges {\rm(}to a random right infinite word in
generators $x_1,\ldots,x_n$ and their inverses\/{\rm)}.
\end{proclaim}

Below, we calculate the reduced words $A_4(\gamma_k)x_4A_4^{-1}(\gamma_k)$ for
a part of the above-considered path in $P_4$.

\begin{equation*}
\begin{array}{lrl}
\gamma_0(x_4)=&
&\!\!\!\!\!x_4\\ \gamma_1(x_4)=&
&\!\!\!\!\!x_4\\ \gamma_2(x_4)=&   x_4^{-1}x_2^{-1}x_1^{-1}x_2
&\!\!\!\!\!x_4x_2^{-1}x_1x_2x_4   \\ \gamma_3(x_4)=&
x_2^{-1}x_1x_2x_3x_2^{-1}x_1^{-1}x_2x_4^{-1}x_2^{-1}x_1^{-1}x_2
&\!\!\!\!\!x_4x_2^{-1}x_1x_2x_4x_2^{-1}x_1x_2x_3^{-1}x_2^{-1}x_1^{-1}x_2   \\
\gamma_4(x_4)=&
\gamma_3(s_{32}^{-1}(x_4))=\gamma_3(x_4)\qquad\qquad\qquad\quad  &  \\
\gamma_5(x_4)=&   \cdots
x_2^{-1}x_1x_2x_3x_2^{-1}x_1^{-1}x_2x_4^{-1}x_2^{-1}x_1^{-1}x_2
&\!\!\!\!\!x_4x_2^{-1}x_1x_2x_4x_2^{-1}x_1x_2x_3^{-1}x_2^{-1}x_1^{-1}x_2
\cdots \\ \gamma_6(x_4)=&   \cdots
x_2^{-1}x_1x_2x_3x_2^{-1}x_1^{-1}x_2x_4^{-1}x_2^{-1}x_1^{-1}x_2
&\!\!\!\!\!x_4x_2^{-1}x_1x_2x_4x_2^{-1}x_1x_2x_3^{-1}x_2^{-1}x_1^{-1}x_2
\cdots
\end{array}
\end{equation*}

There are several other ways to observe the stabilization via Artin's
representation. For example, consider the ratio
$$\frac{\#_{x_p}A_i(\gamma)}{\#_{x_q}A_i(\gamma)},$$ where
$\#_{x_r}A_i(\gamma)$ denotes the number of occurrence of the symbol $x_r$ in
the reduced word $A_i(\gamma)$. It turns out that this ratio is a stable
functional on the braid group.

Note also that for a path $\gamma_1,\gamma_2,\ldots$ the corresponding
sequence of elements $\gamma_k(x_i)$ in $F_n$ is not a.s. convergent (in the
usual sense of the hyperbolic compactification). However, let us consider the
factor-space of $F_n$ by the left translation by the element $x_1x_2\cdots
x_n$. It turns out that in this factor-space the image of $\gamma_k(x_i)$ is
an a.s. convergent sequence.

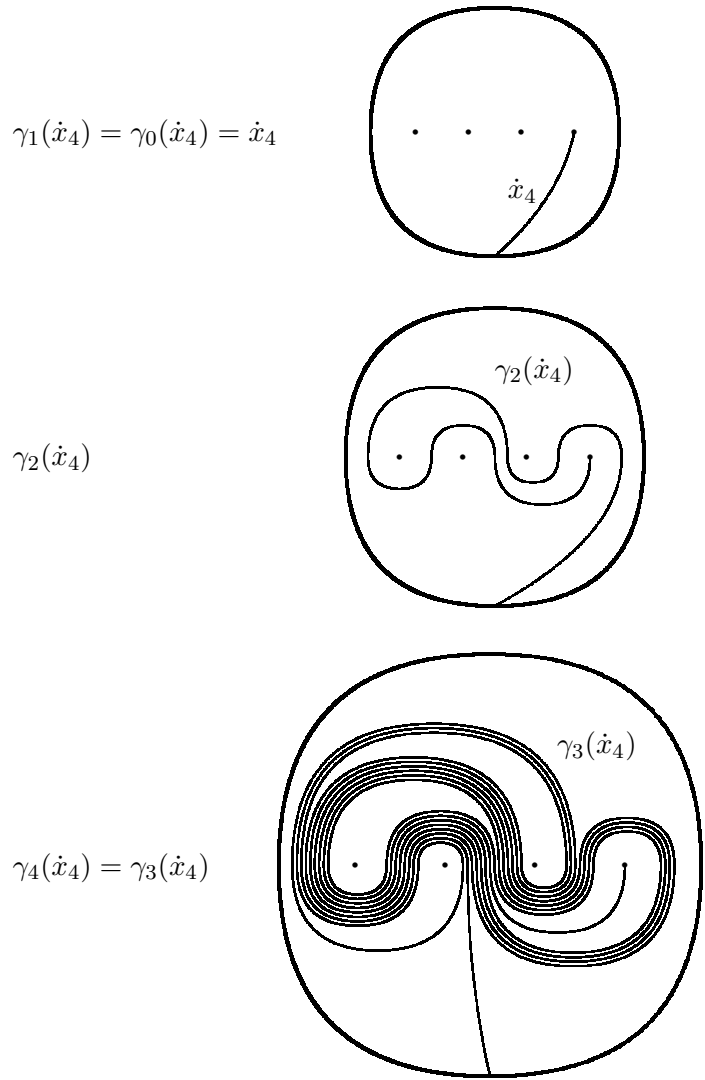
\begin{figure}[ht]
\unitlength=0.1pt
\centerline {\begin{picture}(1000,1100)  
\linethickness{1pt}
\put(500,500)
{
\bezier500(-470,0)(-470,470)(0,470)\bezier200(0,470)(470,470)(470,0)
\bezier500(-470,0)(-470,-470)(0,-470)\bezier200(0,-470)(470,-470)(470,0)
}
\linethickness{0.5pt}
\put(0,500)
{
\put(200,0){\circle*{20}}\put(400,0){\circle*{20}}\put(600,0){\circle*{20}}\put(800,0){\circle*{20}}
\bezier200(500,-470)(750,-250)(800,0)
\put(550,-250){$\dot{x}_4$}
}
\end{picture}}
\unitlength=0.12pt
\centerline { \begin{picture}(1000,1100)  
\linethickness{1pt}
\put(500,500)
{
\bezier500(-470,0)(-470,470)(0,470)\bezier200(0,470)(470,470)(470,0)
\bezier500(-470,0)(-470,-470)(0,-470)\bezier200(0,-470)(470,-470)(470,0)
}
\linethickness{0.5pt}
\put(0,500)
{
\put(200,0){\circle*{20}}\put(400,0){\circle*{20}}\put(600,0){\circle*{20}}\put(800,0){\circle*{20}}
\put(650,0){\bezier200(-150,0)(-150,-150)(0,-150)\bezier200(0,-150)(150,-150)(150,0)}
\put(400,0){\bezier200(-100,0)(-100,100)(0,100)\bezier200(0,100)(100,100)(100,0)}
\put(200,0){\bezier200(-100,0)(-100,-100)(0,-100)\bezier200(0,-100)(100,-100)(100,0)}
\put(320,0){\bezier200(-220,0)(-220,220)(0,220)\bezier200(0,220)(220,220)(220,0)}
\put(620,0){\bezier200(-80,0)(-80,-80)(0,-80)\bezier200(0,-80)(80,-80)(80,0)}
\put(800,0){\bezier200(-100,0)(-100,100)(0,100)\bezier200(0,100)(100,100)(100,0)}
\bezier200(500,-470)(900,-250)(900,0)
\put(500,250){$\gamma_2(\dot{x}_4)$}
}
\end{picture} }
\unitlength=0.17pt
\centerline{\begin{picture}(1000,1050)  
\linethickness{1pt}
\put(500,500)
{
\bezier500(-470,0)(-470,470)(0,470)\bezier200(0,470)(470,470)(470,0)
\bezier500(-470,0)(-470,-470)(0,-470)\bezier200(0,-470)(470,-470)(470,0)
}
\linethickness{0.5pt}
\put(0,500)
{
\put(200,0){\circle*{9}}\put(400,0){\circle*{9}}\put(600,0){\circle*{9}}\put(800,0){\circle*{9}}
\put(375,0)
{\bezier200(-315,0)(-315,315)(0,315)\bezier200(0,315)(315,315)(315,0)
\bezier200(-305,0)(-305,305)(0,305)\bezier200(0,305)(305,305)(305,0)
\bezier200(-295,0)(-295,295)(0,295)\bezier200(0,295)(295,295)(295,0)}
\put(330,0)
{\bezier200(-240,0)(-240,240)(0,240)\bezier200(0,240)(240,240)(240,0)
\bezier200(-230,0)(-230,230)(0,230)\bezier200(0,230)(230,230)(230,0)
\bezier200(-220,0)(-220,220)(0,220)\bezier200(0,220)(220,220)(220,0)
\bezier200(-210,0)(-210,210)(0,210)\bezier200(0,210)(210,210)(210,0)
\bezier200(-200,0)(-200,200)(0,200)\bezier200(0,200)(200,200)(200,0)
\bezier200(-190,0)(-190,190)(0,190)\bezier200(0,190)(190,190)(190,0)}
\put(390,0)
{\bezier200(-120,0)(-120,120)(0,120)\bezier200(0,120)(120,120)(120,0)
\bezier200(-110,0)(-110,110)(0,110)\bezier200(0,110)(110,110)(110,0)
\bezier200(-100,0)(-100,100)(0,100)\bezier200(0,100)(100,100)(100,0)
\bezier200(-090,0)(-090,090)(0,090)\bezier200(0,090)(090,090)(090,0)
\bezier200(-080,0)(-080,080)(0,080)\bezier200(0,080)(080,080)(080,0)
\bezier200(-070,0)(-070,070)(0,070)\bezier200(0,070)(070,070)(070,0)
\bezier200(-060,0)(-060,060)(0,060)\bezier200(0,060)(060,060)(060,0)
\bezier200(-050,0)(-050,050)(0,050)\bezier200(0,050)(050,050)(050,0)}
\put(805,0)
{\bezier200(-105,0)(-105,105)(0,105)\bezier200(0,105)(105,105)(105,0)
\bezier200(-095,0)(-095,095)(0,095)\bezier200(0,095)(095,095)(095,0)
\bezier200(-085,0)(-085,085)(0,085)\bezier200(0,085)(085,085)(085,0)
\bezier200(-075,0)(-075,075)(0,075)\bezier200(0,075)(075,075)(075,0)}
\put(250,0)
{\bezier200(-190,0)(-190,-190)(0,-190)\bezier200(0,-190)(190,-190)(190,0)}
\put(205,0)
{\bezier200(-135,0)(-135,-135)(0,-135)\bezier200(0,-135)(135,-135)(135,0)
\bezier200(-125,0)(-125,-125)(0,-125)\bezier200(0,-125)(125,-125)(125,0)
\bezier200(-115,0)(-115,-115)(0,-115)\bezier200(0,-115)(115,-115)(115,0)
\bezier200(-105,0)(-105,-105)(0,-105)\bezier200(0,-105)(105,-105)(105,0)
\bezier200(-095,0)(-095,-095)(0,-095)\bezier200(0,-095)(095,-095)(095,0)
\bezier200(-085,0)(-085,-085)(0,-085)\bezier200(0,-085)(085,-085)(085,0)
\bezier200(-075,0)(-075,-075)(0,-075)\bezier200(0,-075)(075,-075)(075,0)
\bezier200(-065,0)(-065,-065)(0,-065)\bezier200(0,-065)(065,-065)(065,0)}
\put(620,0)
{\bezier200(-050,0)(-050,-050)(0,-050)\bezier200(0,-050)(050,-050)(050,0)
\bezier200(-060,0)(-060,-060)(0,-060)\bezier200(0,-060)(060,-060)(060,0)
\bezier200(-070,0)(-070,-070)(0,-070)\bezier200(0,-070)(070,-070)(070,0)
\bezier200(-080,0)(-080,-080)(0,-080)\bezier200(0,-080)(080,-080)(080,0)
\bezier200(-090,0)(-090,-090)(0,-090)\bezier200(0,-090)(090,-090)(090,0)
\bezier200(-100,0)(-100,-100)(0,-100)\bezier200(0,-100)(100,-100)(100,0)
\bezier200(-110,0)(-110,-110)(0,-110)\bezier200(0,-110)(110,-110)(110,0)}
\put(685,0)
{\bezier200(-225,0)(-225,-225)(0,-225)\bezier200(0,-225)(225,-225)(225,0)
\bezier200(-215,0)(-215,-215)(0,-215)\bezier200(0,-215)(215,-215)(215,0)
\bezier200(-205,0)(-205,-205)(0,-205)\bezier200(0,-205)(205,-205)(205,0)
\bezier200(-195,0)(-195,-195)(0,-195)\bezier200(0,-195)(195,-195)(195,0)}
\put(650,0){\bezier200(-150,0)(-150,-150)(0,-150)\bezier200(0,-150)(150,-150)(150,0)}
\bezier200(500,-470)(450,-250)(450,0)
\put(650,250){$\gamma_3(\dot{x}_4)$}
} 
\end{picture} }
\unitlength=1pt
\centerline{\begin{picture}(1,1)
\put(-180,373){$\gamma_1(\dot{x}_4)=\gamma_0(\dot{x}_4)=\dot{x}_4$}
\put(-180,250){$\gamma_2(\dot{x}_4)$}
\put(-180,95){$\gamma_4(\dot{x}_4)=\gamma_3(\dot{x}_4)$}
\end{picture} }
\caption{The action of $4$-braids on the geodesic $\dot{x}_4$.}
\label{disc}
\end{figure}

\paragraph{The braid group as the mapping class group.}
Next, we consider the well-known interpretation of the braid group in the form
of the mapping class group of the punctured disc. Under this interpretation,
the action of the braid group on the fundamental group of the punctured disc
is exactly the Artin representation. It is convenient to regard the elements
of the type $\gamma(x_i)=A_i(\gamma)x_iA_i^{-1}(\gamma)$ in this
free--fundamental group as the homotopy classes of curves that come out of the
$i$-th puncture and end in a base point on the boundary. It is also convenient
to endow our punctured disc with a hyperbolic metric. Then each of these
homotopy classes is represented by a unique geodesic (see Figure~\ref{disc}).
The braid group acts on the set of such geodesics. In this construction, the
described above (left) stabilization of the sequence $A_n^{-1}(\gamma_k)$ for
a path $\gamma_1,\gamma_2,\ldots$ transforms into the stabilization of the
sequence of geodesics (here, a stable functional is the behaviour of a
geodesic after it comes out of the puncture (see Figure~\ref{disc}).

We remark that the study of the interpretation of the braid group in the form
of the mapping class group leads finally to a $\rm{mod}\,0$-isomorphism between
the PF-boundary of the braid group described in this paper (the boundary of the
free group) and the PF-boundary of the braid group in terms of~\cite{KM,FM}
(the Thurston boundary of Teichm\"uller space of the punctured sphere).

\paragraph{Representations of braid groups in the group of homeomorphisms of the circle.}

With the aid of the above-mentioned interpretation of the braid group in the
form of the mapping class group, we can construct representations of the braid
group in the group of homeomorphisms of the circle. Consider, for example, the
set of all oriented geodesics (on the punctured disc with a hyperbolic metric)
that come out of a certain puncture. This set is naturally endowed with the
topology of the circle. The pure braid group acts on this ``circular'' space
of geodesics. We have thus converted the pure braid group into a group, which
acts on the circle. For this representation, our results yield corresponding
corollaries. In particular, we see that the circle is a $\mu$-boundary of the
pure braid group.

Using the punctured disc, we also can let the whole braid group to act on the
circle. For example, we can define an action of the braid group on the set of
all oriented geodesics that start from the boundary at right angle to it.
(Obviously, this set is endowed with the topology of the circle.) For this
action, Theorems~1 and~2 also produce certain corollaries.

\subsection{Other normal forms of braids}
The Markov--Ivanovsky normal form is not a unique stable form for braids. In
fact, the results of this paper imply the stability for several other normal
forms, which are essentially different from that of Markov--Ivanovsky. At the
other hand, for some well-known normal forms (of braids) the problem of
stability is still open.

Going into the details, first we mention a stable normal form described
in~\cite{Mal1}. Actually, this form belongs to an extensive collection of
stable normal forms of geometrical nature. The stability of this form follows
from the stability of the Markov--Ivanovsky form. In~\cite{Bre}, Bressaud
presented a new normal form for braids, which is to a certain extent cognate to
the form of~\cite{Mal1}. The Bressaud form is not stable with respect to the
right-hand random walk. Nevertheless, this form can be easily converted into a
stable one by shifting its leftmost part (which relates with the center of the
braid group) to the right. The stability of this converted Bressaud form can
also be deduced from the results of this work. There were attempts to find a
stable form with all words $\si$-positive or $\si$-negative in Dehornoy's sense
(see~\cite{DDRW} for definitions). It turns out that such a form does not exist
(for the random walk with a symmetric distribution). At the same time, we can
construct a stable normal form composed of two parts, one of which is a
$\si$-positive or a $\si$-negative word while the second is in a certain sense
elementary. For the Garside normal form, Thurston's normal form~\cite{Eps},
Birman-Ko-Lee normal form~\cite{BKL}, and some other known forms
(see~\cite{Deho} for examples) the stability problem\footnote{It is obvious
that some of the aforementioned forms are not stable with their customary
definitions. Here, by the stability problem for a normal form we mean whether a
normal form or some of its natural simple transformations is stable. Examples
of such transformations (which usually include permutations of the parts of a
form under investigation) are mentioned in the above remarks on the Bressaud
normal form stability and also in the remarks on the Markov--Ivanovsky form
definition in Section~\mmBG\ below.} is open.

\section{The braid group and the Markov--Ivanovsky normal form. Definitions}

In his seminal work~\cite{Mar}, A.A.\,Markov presented a normal form for the
elements of the braid group. This form is based on a normal series of the pure
braid group, which consists of free groups with decreasing ranks. The first
free group in this series is a {\it special\/} free subgroup of the braid
group. The pure braid group $P_n$ is a semidirect product of the special free
subgroup by a subgroup isomorphic to the pure braid group $P_{n-1}$. By
induction we therefore obtain a normal form for the pure braid group. This form
is composed by a sequence of elements in the free subgroups of the
aforementioned normal series. A.A.\,Markov annotates the theorems about this
normal form and about the normal series with the last name of A.\,Ivanovsky.
Unfortunately, neither the work~\cite{Mar} nor other papers of A.A. contain any
detailed references on A.\,Ivanovsky\footnote{A.\,Ivanovsky was a post-graduate
student of A.A.\,Markov. He was killed during World War II.}. In later works on
the braid groups this normal series was rediscovered but it seems that the
first appearance of the form was in~\cite{Mar}, so we prefer the term ``the
Markov--Ivanovsky normal form''.

The {\it Artin braid group\/} of rank $n$ is determined by the presentation
\begin{equation*}
B_n := \langle \sigma_1, \dots, \sigma_{n-1} \mid \sigma_i \sigma_j = \sigma_j
\sigma_i, \ |i-j| \ge 2; \quad \sigma_i \sigma_{i+1} \sigma_i =\sigma_{i+1}
\sigma_i \sigma_{i+1} \rangle.
\end{equation*}
In $B_n$, we consider the set of elements $\{s_{ji}, 1\le i < j \le n\}$, where
$$s_{ji}:=\si_{j-1}\si_{j-2}\cdots\si_{i+1}\si_i^2\si_{i+1}^{-1}\cdots\si_{j-2}^{-1}\si_{j-1}^{-1}.$$
Suppose $m$ is in $\{2,\ldots n\}$; then the set $\{s_{ji}, 1\le i < j \le m\}$
generates the {\it pure braid group\/} $P_m$ in $B_n$. The collection
$\{s_{mi}, 1\le i < m\}$ generates a free subgroup of rank $m-1$, which we
denote by $F_{m-1}$. $F_{m-1}$ is a normal subgroup in $P_m$. The subgroups
$P_2\cong \Z$ and $F_1\cong \Z$ coincide. For each $k\in\{3,\ldots n\}$, the
subgroup $P_k$ is a semidirect product of $F_{k-1}$ by $P_{k-1}$: $$
P_k=F_{k-1} \rtimes P_{k-1}. $$

Therefore, we have $$P_n= F_{n-1} \rtimes (F_{n-2}\rtimes (F_{n-3}\rtimes
(F_{n-4}\rtimes\cdots\rtimes{F_1}))).$$

In particular, each element $\gamma$ of the pure braid group $P_n$ is uniquely
presented in the form $$\gamma = \gamma_{n-1}
\gamma_{n-2}\cdots\gamma_2\gamma_1, \quad\text{where}\quad \gamma_i \in F_i.$$
Every element of the free group $F_i$ is represented by a unique {\it
reduced\/} word in generators $\{s_{(i+1)j}, 1\le j <i+1\}$ and their inverses
(the reader may find an explicit definition of reduced words in
Subsection~\mmF.1 below).

The {\it Markov--Ivanovsky normal form in the pure braid group $P_n$\/} is the
mapping $$ \mathfrak{I}_P:\gamma \mapsto V_{n-1}\cdots V_1, $$ where $V_i$ is
the reduced word for the element $\gamma_i \in F_i$.

In order to define the Markov--Ivanovsky normal form in the braid group $B_n$,
we recall that $P_n$ is a normal subgroup of index $n!$ in $B_n$. Let $\Pi_n
\subset B_n$ be a set of coset representatives for $P_n$ in $B_n$. Then each
element $\beta \in B_n$ can be uniquely written in the form $\gamma_\beta
\pi_{\beta}$, where $\gamma_\beta\in P_n$ and $\pi_{\beta} \in \Pi_n$. The {\it
Markov--Ivanovsky normal form in the braid group $B_n$\/} is the
mapping\footnote{In~\cite{Mar}, Markov fixes a certain specific set of coset
representatives (for $P_n$ in $B_n$) to define the normal form. In our case,
the choice of such a set is not of importance. Furthermore, in~\cite{Mar} the
sequence order for ``components'' of the normal form differs from the sequence
order used above: Markov wrote a braid $\beta \in B_n$ in the form $$
\mathfrak{I}'_B(\beta) = \beta'_1\beta'_2\cdots\beta'_{n-2}\beta'_{n-1}
\pi'_{\beta}, \quad\text{where}\quad \beta'_i \in F_i, \quad \pi'_{\beta} \in
\Pi'_n.$$ Clearly, in certain cases it is not important which sequence order we
choose for the definition of the normal form. Apparently, in~\cite{Mar} the
sequence order was chosen arbitrarily. We observe that in~\cite{Vershinin}, for
example, the Markov--Ivanovsky normal form is described with the same sequence
order as ours. In the context of the stability, the choice of a sequence order
does matter: the normal form $\mathfrak{I}_B$ is stable with respect to the
right random walk, while the normal form $\mathfrak{I}'_B$ is not stable (for
both left and right random walks).}
$$\mathfrak{I}_B:\beta\mapsto\mathfrak{I}_P(\gamma_\beta)\pi_{\beta}.$$

\section{Random walks on groups. Definitions}

\paragraph{\mmRW.1. Random walks on groups.}
Let $G$ be a countable group, and $\mu$ be a probability measure on~$G$. We say
that $\mu$ is {\it nondegenerate\/} if its support generates $G$ as a
semigroup. The (right-hand) {\it random walk on $G$ determined by the
distribution $\mu$\/}, or the {\it random $\mu$-walk\/}, is the time
homogeneous Markov chain with the state space $G$, with the transition
probabilities $P(g,h) = \mu(g^{-1}h)$, and with the initial distribution
concentrated at the group identity $e$. A realization of this process is called
a {\it path\/} of the random walk. If $\mu$ is a probability measure on $G$, we
denote by $\mi$ the associated Markov measure on the path space $G^{\Z_+}$.

\paragraph{\mmRW.2. Stable normal forms~\cite{Ver}.}
Let $S$ be a subset of $G$. Suppose that $S$ generates $G$ as a semigroup. We
denote by $S^*$ the set of all finite {\it words\/} over the {\it alphabet\/}
$S$. A {\it normal form\/} in the group $G$ is a mapping $\mathfrak{N}: G \to
S^*$ such that $\operatorname{pr} \circ \mathfrak{N} = \id_G$, where
$\operatorname{pr}$ is the natural projection from $S^*$ onto $G$, that is, the
projection that takes the word $w_1 w_2 \dots w_k\in S^*$ to the element
$w_1\cdot w_2\dots w_k\in G$. We say that a sequence $\{V_i\}_{i\in\Z_+}$ of
words in $S^*$ {\it converges at infinity\/} if for every $k\in\N$ there exists
$N\in\N$ such that for each $j> N$ the length of the word $V_j$ is greater then
$k$ and, moreover, the initial subwords of length $k$ (the $k$-{\it
prefixes\/}) of words $V_j$, $V_N$ coincide. A normal form $\mathfrak{N}:G \to
S^*$ is said to be {\it stable with respect to the random $\mu$-walk\/}
($\mu$-{\it stable\/}) if for $\mi$-a.e. path $\tau=\{\tau_i\}_{i\in\Z_+}$ the
sequence of words $\{\mathfrak{N}(\tau_i)\}_{i\in\Z_+}$ converges at
infinity\footnote{It is interesting to investigate also the notion of the {\it
weak stability\/} for normal forms. (To obtain the definition of the weak
stability, use the convergence in measure instead of the a.e. convergence.)}. A
normal form is said to be {\it stable\/} if it is $\mu$-stable for every
nondegenerate measure $\mu$.

\paragraph{\it Remark.}
It is a too restrictive approach, if we define normal forms as words, and
stabilization as the convergence of words in the word metric. A more efficient
approach to the stable normal forms can be obtained by generalizing the notion
of ``word''. A usual (linear) word of length $n$ is a mapping of $\{1, \dots
n\}\subset \N$ to an alphabet $S$. Yet we can consider the union of several
copies of $\N$ (the ``vector words''). Such a generalization of the notion of
normal forms is necessary if our group $G$ is the product of several copies of
free groups or braid groups and so on. Accordingly, under this approach the
limits of a stable normal form are collections of several infinite linear
words.

However, this interpretation must be further generalized. For example, in the
case of the free meta-abelian group (see~\cite{Ver}), it is convenient to
define a normal form of an element as a configuration of paths of a certain
shape on the lattice with the operation of concatenation. Here, the limits are
the infinite configurations. This interpretation for normal forms is especially
useful if our group is realized as the fundamental group of a certain complex.
This leads us to a notion of {\it geometric boundary of the random walk on a
group}.

\paragraph{\mmRW.3. $\mu$-boundaries.}
Recall that an {\it action\/} of a group $G$ on a topological
space $M$ is a homomorphism from $G$ to the group
$\operatorname{Homeo}(M)$ of all homeomorphisms of $M$.
A space endowed with an action of a group $G$ is called {\it $G$-space.} We
denote by $\PP(G)$ the space of all probability measures on $G$. The space of
all regular Borel probability measures on a topological space $M$ is denoted by
$\PP(M)$. $\PP(M)$ is endowed with the weak$^*$ topology (with respect to the
bounded continuous functions on $M$). A measure concentrated at a point is
called a {\it point measure\/}. Recall also that an action of $G$ on $M$
determines an action of $G$ on $\PP(M)$, $g\theta(E)=\theta (g^{-1}E)$.

Let $\mu\in\PP(G)$. A measure $\nu\in\PP(M)$ is said to be {\it
$\mu$-stationary\/} if $$ \mu * \nu = \int_G g\nu\,d\mu(g)=\nu. $$

Let $M$ be a $G$-space, let $\mu\in\PP(G)$, and let $\nu\in\PP(M)$ be a
$\mu$-stationary measure. A~pare $(M, \nu)$ is called a {\it $\mu$-boundary\/}
for $G$ if for a.e. path $\tau=\{\tau_i\}_{i\in\Z_+}$ of the random $\mu$-walk
the sequence of measures $\{\tau_i (\nu)\}_{i\in\Z_+}$ converges to a point
measure $\delta_{w(\tau)}$, where $w(\tau)\in M$. A $G$-space $M$ is said to be
{\it $\mu$-proximal\/} if for every $\mu$-stationary measure $\nu\in\PP(M)$ the
pair $(M, \nu)$ is a $\mu$-boundary. A $G$-space $M$ is said to be {\it
mean-proximal\/} if it is $\mu$-proximal for every nondegenerate measure $\mu$.

We say that a $\mu$-boundary $(M, \nu)$ of a pair $(G,\mu)$ is a {\it
Poisson--Furstenberg boundary\/} (a {\it PF-boundary\/}) if $(M, \nu)$ is a
{\it maximal\/} $\mu$-boundary, i.e., if each $\mu$-boundary $(M_1, \nu_1)$ of
the pair $(G,\mu)$ is a factor-space of $(M,\nu)$ (as a measure space with a
group action, disregarding the topology).

\section{Random walks on groups. Lemmas}

This section contains several general statements related to random walks on
groups.

\begin{proclaim}{\mmRWL.1. Proposition \cite{Fu3}} 
Assume that a countable group $G$ acts on a compact metric space $M$. Then for
every measure $\mu\in\PP(G)$ the set of $\mu$-stationary measures in $\PP(M)$
is nonempty. \qed
\end{proclaim}

\begin{proclaim}{\mmRWL.2. Theorem \cite{Fu3}} 
Assume that a countable group $G$ acts on a compact metric space $M$. Let
$\nu\in\PP(M)$ be a $\mu$-stationary measure for a measure $\mu\in\PP(G)$.
Then for a.e. path $\tau=\{\tau_i\}_{i\in\Z_+}$ of the random $\mu$-walk the
sequence $\{\tau_i (\nu)\}_{i\in\Z_+}$ converges to a measure
$\lambda(\tau)\in\PP(M)$, and $$ \int \lambda(\tau) d\mi(\tau)=\nu. \qed $$
\end{proclaim}

\begin{proclaim}{\mmRWL.3. Lemma}
Assume that a countable group $G$ acts on a compact metric space $M$. Let
$\mu\in\PP(G)$. Assume that $M$ is $\mu$-proximal. Then a $\mu$-stationary
measure on $M$ is unique.
\end{proclaim}

\begin{proof}
Let $\nu_1$, $\nu_2$ be two $\mu$-stationary measures on $M$. Then, by the
definition of $\mu$-stationary measures, the measure $\nu_0:=\frac{1}{2}\nu_1 +
\frac{1}{2}\nu_2$ is also $\mu$-stationary. Since $M$ is $\mu$-proximal, it
follows that for a.e. path $\tau:=\{\tau_i\}_{i\in\Z_+}$ of the random
$\mu$-walk the corresponding sequence
$\{\tau_i\nu_0\}_{i\in\Z_+}=\{\frac{\tau_i\nu_1}{2} +
\frac{\tau_i\nu_2}{2}\}_{i\in\Z_+}$ of measures converges to a point measure
$\delta_{w}$, where $w=w(\tau)\in M$. This implies that the sequences
$\{\frac{\tau_i\nu_1}{2}\}_{i\in\Z_+}$ and
$\{\frac{\tau_i\nu_2}{2}\}_{i\in\Z_+}$ converge to the measure
$\frac{\delta_{w}}{2}$. Then by Theorem~\mmRWL.2 we have $$ \nu_0=\int
\delta_{w(\tau)} d \mi(\tau)=\nu_1=\nu_2. $$
\end{proof}

\begin{proclaim}{\mmRWL.4. Lemma \cite{KM}} 
Assume that a countable group $G$ acts on a space $M$. Let $\mu$ be a
nondegenerate measure on $G$, and $\nu\in\PP(M)$ be a $\mu$-stationary measure.
Suppose $E$ is a measurable subset of $M$ such that for every $g\in G$ we have
either $g(E)=E$ or $g(E)\cap E =\varnothing$, and, moreover, there are an
infinite number of pairwise disjoint sets of the form $g(E)$, where $g\in G$.
Then $\nu(E)=0$.

In particular, if for each point $x\in M$ the orbit $G(x)$ is infinite, then
$\nu$ is continuous.
\end{proclaim}

\begin{proclaim}{\mmRWL.5. Lemma (on absolute continuity)}
Assume that a countable group $G$ acts on a space $M$. Let $\mu$ be a
nondegenerate measure on $G$, and let $\nu\in\PP(M)$ be a $\mu$-stationary
measure. Then for each $g\in G$ the measure $g\nu$ is absolutely continuous
with respect to the measure $\nu$. Moreover, we have

\begin{itemize}

\item[{\rm 1)}] For every $g\in G$ there exists a constant $C'_g$ such that for each
measurable subset $E\subset M$ we have $$ C'_g \nu(E)\ge g\nu(E). $$

\item[{\rm 2)}] For every $g\in G$ there exists a constant $C''_g >0$ such that for each
measurable subset $E\subset M$ with $E \cup gE = M$ we have $$ \nu(E)\ge
C''_g. $$

\end{itemize}
\end{proclaim}

\begin{proof}
\begin{itemize}

\item[1)] Since $\mu$ is nondegenerate, it follows that for an element $g\in G$ there
is a number $s\in\N$ such that $\mu_s(g)>0$, where $\mu_s$ denotes the
$s$-fold convolution of $\mu$. Since $\nu$ is $\mu$-stationary, it follows
that $$ \nu(E)=\sum_{h\in G}h\nu(E)\mu_s(h)\ge g\nu(E)\mu_s(g). $$ Now, we set
$C'_g:=1/\mu_s(g)$.

\item[2)] The condition $E \cup gE = M$ implies that $g^{-1}E \cup E = g^{-1}M=M$.
From this we obtain $$ \nu(E)+g\nu(E)=\nu(E)+\nu(g^{-1}E)\ge \nu(E \cup
g^{-1}E)=\nu(M)=1. $$ By assertion 1) we have $C'_g \nu(E)\ge g\nu(E)$.
Consequently, $$ (1+C'_g)\nu(E)\ge \nu(E)+g\nu(E)\ge 1. $$ To complete the
proof, we set $C''_g:=\frac{1}{1+C'_g}$.

\end{itemize}
\end{proof}

\paragraph{Definition \rm{(cofinite subsets)}.}
We say that a subset $H$ of a group $G$ is {\it cofinite\/} if there exists a
finite set $J\subset G$ such that $HJ=G$. The following result immediately
follows from the standard facts of martingale theory.

\begin{proclaim}{\mmRWL.6. Lemma (on cofinite subsets)}
Let $G$ be a countable group, $H\subset G$ be a cofinite subset, $\mu\in\PP(G)$
be a nondegenerate measure. Then a.e. path of the random $\mu$-walk hits the
subset $H$ an infinite number of times. Furthermore, for each element $g_0\in
G$ there exists $k\in\N$ such that for a.e. path $\tau=\{\tau_i\}_{i\in\Z_+}$
of the random $\mu$-walk the set $$ L(\tau):=\{l\in\N\ |\ \tau_l\in H,\
\tau_{l+k}=\tau_l g_0\} $$ is infinite.\qed
\end{proclaim}

\section{A sufficient condition of $\mu$-proximality}

In this section, we describe a sufficient condition for a metric $G$-space to
be $\mu$-proximal. This condition will be used in the proof of Theorem~1
below.

\paragraph{\mmPRIZ.1. $\varepsilon$-contracting collections of elements.}
Let $M$ be a metric space with a metric~$d$. Suppose that a group $G$ acts on
$M$ by homeomorphisms. An element $g\in G$ is said to be {\it
$\varepsilon$-contracting for a measure $\lambda\in\PP(M)$\/} if there is a
ball $$B_\varepsilon(x):=\{y\in M|\ d(x,y)\le \varepsilon\}$$ such that $$
g\lambda\big(B_\varepsilon(x)\big) \ge 1-\varepsilon. $$

We say that a collection $K_1 \subset G$ is {\it totally
$\varepsilon$-contracting\/} (for a metric $G$-space $(M,d)$) if for each
measure $\lambda\in\PP(M)$ the collection $K_1$ contains an
$\varepsilon$-contracting element for~$\lambda$.

A collection $K_2\subset G$ is said to be {\it universally\/}
$\varepsilon$-{\it contracting\/} (for $(M,d)$) if for each $h\in G$ the
collection $hK_2$ is totally $\varepsilon$-contracting.

\paragraph{\mmPRIZ.2. Remark.}
If a collection $K_1\subset G$ is totally $\varepsilon$-contracting, then for
each $h\in G$ the collection $K_1h$ is also totally $\varepsilon$-contracting.
To prove this, let $\lambda\in\PP(M)$ be an arbitrary measure. Then $K_1$
contains an $\varepsilon$-contracting element $g$ for the measure $h\lambda$.
This means that the element $gh\in K_1h$ is $\varepsilon$-contracting for
$\lambda$.

\paragraph{\mmPRIZ.3. Remark.}
The previous remark immediately implies that a set $K_2\subset G$ is
universally $\varepsilon$-contracting if and only if for each $h\in G$ the set
$hK_2h^{-1}$ is totally $\varepsilon$-contracting.

\begin{proclaim}{\mmPRIZ.4. Proposition (a suffitient condition of $\mu$-proximality)}
Assume that a countable group $G$ acts on a compact metric space $M$. Suppose
that for an arbitrary small $\varepsilon>0$ the group $G$ contains a finite
universally $\varepsilon$-contracting collection for $M$. Then for each
nondegenerate measure $\mu$ on $G$, $M$ is $\mu$-proximal.
\end{proclaim}

The proof of Proposition~\mmPRIZ.4 is based upon the following lemma.

\begin{proclaim}{\mmPRIZ.5. Lemma}
Assume that a group $G$ acts on a metric space $M$, $\varepsilon>0$, and $G$
contains a finite universally $\varepsilon$-contracting collection $K$ for~$M$.
Then for an arbitrary measure $\lambda\in\PP(M)$ the set
$H:=H(\varepsilon,\lambda)$ of all $\varepsilon$-contracting elements for
$\lambda$ is cofinite in~$G$.
\end{proclaim}

\begin{proof}[Proof of the lemma]
Let $g\in G$ be an arbitrary element. Then, by the definition of a universally
$\varepsilon$-contracting collection, the collection $gK$ is totally
$\varepsilon$-contracting. In other words, the sets $gK$ and $H$ have a
nonempty intersection, whence $g\in HK^{-1}$. Then, since $g$ is arbitrary, it
follows that $G=HK^{-1}$. Therefore, since the set $K^{-1}$ is finite, it
follows that $H$ is, indeed, cofinite.
\end{proof}

\begin{proof}[Proof of Proposition~{\rm\mmPRIZ.4}]
Let $\mu\in\PP(G)$ be an arbitrary nondegenerate measure, and let $\nu$ be a
$\mu$-stationary measure in $\PP(M)$. By Lemma~\mmPRIZ.5, for each
$\varepsilon>0$ the set of all $\varepsilon$-contracting elements for $\nu$ is
cofinite in~$G$. Then Lemma~\mmRWL.6 implies that for a.e. path
$\tau=\{\tau_i\}_{i\in\Z_+}$ of the random $\mu$-walk there exists a sequence
$\{i_k\}_{k\in\Z_+}\subset \N$ such that for each $k\in\N$ the element
$\tau_{i_k}$ is $1/k$-contracting for the measure $\nu$. Since $M$ is compact,
it follows that the sequence $\{\tau_{i_k}(\nu)\}_{k\in\Z_+}$ has a subsequence
that converges to a point measure. This means that for a.e. path $\tau$ the
sequence $\{\tau_i (\nu)\}_{i\in\Z_+}$ has a subsequence converging to a point
measure. At the same time, by Theorem~\mmRWL.2 the sequence $\tau_i(\nu)$ a.s.
converges to a certain measure. Consequently, $\{\tau_i (\nu)\}_{i\in\Z_+}$
a.s. converges to a point measure, i.e. the pair $(M, \nu)$ is a $\mu$-boundary
if $(G, \mu)$. Therefore, since $\nu$ is an arbitrary $\mu$-stationary measure,
it follows by definition that $M$ is $\mu$-proximal, as required.
\end{proof}

\section{The free group. Definitions}

This section contains several definitions and well-known facts about the free
group, its boundary, and the action of the free group on its boundary. We use
standard terminologies of combinatorial group theory and hyperbolic group
theory (see~\cite{Gr,LSh}).

\paragraph{\mmF.1. The free group and words.}
Let $F$ be a free group of rank $n\ge 2$ with generators $u_1, \ldots, u_n$.
We will denote the set of these generators and their inverses by $U$. Let
$U^*$ be the set of all finite words over the alphabet $U$, and let
$U^{\infty}$ denote the set of all right infinite words over the same
alphabet.

If $W=w_1w_2\ldots w_k\in U^*$, then we denote by $W^{-1}$ the word
$w_k^{-1}\ldots w_2^{-1}w_1^{-1}\in U^*$. If $r\in\{1,\ldots,k\}$, the word
$w_1\ldots w_r$ is called the {\it initial subword\/} of length $r$, or the
$r$-{\it prefix\/}, of the word $W$, while the word $w_{k-r+1}\ldots w_k$ is
the {\it terminal subword\/} of length $r$, or the $r$-{\it suffix\/}, of~$W$.
For $q\in\Z$, we denote by $W^q$ the $|q|$-fold concatenation of the
word~$W^{\operatorname{sign}(q)}$. By $W^{+\infty}$ and $W^{-\infty}$ we
denote the right infinite words that are the limits of the sequences
$\{W^q\}_{q\in\N}$ and $\{W^{-q}\}_{q\in\N}$, respectively.

A word over the alphabet $U$ is said to be {\it reduced\/} if it contains no
inverse pair subwords. Each element $x\in F$ has a unique reduced
representative in $U^*$. We denote this representative by~$\la x\ra$. Clearly,
the mapping $x \mapsto \la x\ra$ is a normal form. For $x\in F$, we put
$|x|:=|x|_U:=|\la x\ra|$. In the free group $F$, the function
$d(x,y):=|x^{-1}y|$ is a metric. This metric is called the {\it word metric\/}
associated with the system of generators $\{u_1, \ldots, u_n\}$.

A word $V\in U^*$ is said to be {\it cyclically reduced\/} if the word $VV$ is
reduced. It can be easily checked that for an element $a\in F$ there is a
unique pair $(A_a, X_a)$ of cyclically reduced word $A_a\in U^*$ and reduced
word $X_a\in U^*$ such that $\la a\ra=X_a A_a X_a^{-1}$. We call the words
$X_a$ and $A_a$ the {\it wing\/} and the {\it core\/} of the element $a$ and
denote them by $\shl\la a \ra$ and $\osn\la a \ra$, respectively. (We remark
that these definitions depend on the system of generators.) It can be easily
seen that {\it the lengths of cores of conjugate elements are equal\/} (i.e.,
for any $a,b\in F$ we have $|\osn\la a \ra|=|\osn\la bab^{-1} \ra|$).

\paragraph{\mmF.2. The boundary of the free group.}
As is customary, we denote by $\dd F$ the boundary of the free non-abelian
group $F$. It is well known that $\dd F$ is homeomorphic to the Cantor set.
$\dd F$ can be defined as the space of ends for the Cayley graph of $(F,\{u_1,
\ldots, u_n\})$. (Recall that the Cayley graph of the free group with a system
of free generators is a tree.)

There exists a natural one-to-one correspondence between the set of points in
$\dd F$ and the set of all reduced words in $U^{\infty}$. We denote by $\la
w\ra$ the right infinite reduced word that corresponds to a point $w\in\dd F$.
Clearly, the formula $x\mapsto \la x\ra$ determines a bijection from the set
of points of the {\it hyperbolic compactification\/} $F\cup\dd F$ to the set
of all reduced words in $U^*\cup U^{\infty}$.

It is known that for each nontrivial element $a\in F$ the sequence
$\{a^i\}_{i\in\N}$ converges in $F\cup\dd F$ to a boundary point. We denote by
$a^{+\infty}$ and $a^{-\infty}$ the limits of the sequences $\{a^i\}_{i\in\N}$
and $\{a^{-i}\}_{i\in\N}$, respectively. Assume that $\la a \ra=XAX^{-1}$,
where $X$ is reduced and $A$ is cyclically reduced; then, clearly, the words
$XA^{+\infty}$ and $XA^{-\infty}$ represent the points $a^{+\infty}$ and
$a^{-\infty}$, respectively.

\paragraph{\mmF.3. On the action on the boundary.}
For each element $a\in F$, the action of $a$ on $F$ by left translation
($x\mapsto ax$) can be uniquely extended to a homeomorphism of the compact
$F\cup\dd F$. The image of a subset $E\subset\dd F$ under the left translation
by $a$ will be denoted by $aE$. It is known that each automorphism
$\psi\in\operatorname{Aut}(F)$ can be uniquely extended to a homeomorphism
$\ov\psi:F\cup\dd F\to F\cup\dd F$. Let $a\in F$, and let
$\psi_a\in\operatorname{Inn}(F)$ denote the automorphism of conjugation by
$a$, i.e., $\psi_a(x)=axa^{-1}$. Then we readily see that for each $w\in\dd F$
we have $\ov\psi_a(w)=aw$.

\paragraph{\mmF.4. The Gromov product.}
Let $W=w_1w_2\ldots$ and $V=v_1v_2\ldots$ be two distinct words (either finite
of right infinite). We set $(W|V):=r-1$, where $r$ is the smallest positive
integer such that $v_r\neq w_r$ (we mean that the inequality $v_r\neq w_r$
holds true, in particular, if we have either $|V|=r-1$ or $|W|=r-1$). For a
finite word $V$ we set $(V|V)=|V|$.

For a pair of distinct points $w$, $v$ in $F\cup\dd F$ we set $(w|v):=(\la
w\ra|\la v\ra)$. For a point $v\in F$ we set $(v|v)=|v|$. We easily check that
for every $x$, $y$ in $F$ the following equality holds:
$$(x|y)=\frac12\big(|x|+|y|-|x^{-1}y|\big).$$ The value $(x|y)$ is called the
{\it Gromov product\/} of elements $x$ and $y$.

On the space $F\cup\dd F$, we define a metric $\rho$ by setting
$$\rho(w,v):=\big((w|v)+1\big)^{-1}$$ for distinct elements $w$, $v$ in
$F\cup\dd F$. (We set~$\rho(v,v)=0$.) It can be easily checked that $\rho$ is,
indeed, a metric. The topology of this metric coincide with the standard
topology of hyperbolic compactification on $F\cup\dd F$. We remark that for
every $x,y,z\in F\cup\dd F$, in the triple $\rho(x,y)$, $\rho(x,z)$,
$\rho(y,z)$ either all values coincide or the two greatest values coincide.

In the metric space $(\dd F,\rho)$, we will denote by $B_\varepsilon(v)$ the
closed ball $\{w\in \dd F|\ \rho(v,w)\le\varepsilon\}$. It immediately follows
from the definitions that for every $k\in\N$ and $w\in\dd F$, the ball
$B_{1/k}(w)$ is the set of all infinite reduced words with a common
$(k-1)$-prefix. This implies that the intersection of every two balls in $(\dd
F,\rho)$ either is empty or coincides with one of these balls. We note also
that if $v,w\in\dd F$ and $\rho(v,w)\le1/k$, then $B_{1/k}(v)=B_{1/k}(w)$.

\begin{proclaim}{\mmF.5. Lemma}
For every nontrivial element $a\in F$ we have
$$
\rho(a^{+\infty},a^{-\infty})=\frac{1}{|\shl\la a\ra|+1}. \eqno(\mmF.1)\qquad
$$
$$
(a|a^{+\infty})> \frac{|a|}{2}. \eqno(\mmF.2)\qquad
$$
$$ \rho(a,a^{+\infty}) <
\frac{2}{|a|}. \eqno(\mmF.3)\qquad
$$
\end{proclaim}

\begin{proof}
We put $X:=\shl\la a\ra$ and $A:=\osn\la a\ra$. Then $\la a\ra=XAX^{-1}$, $\la
a^{+\infty}\ra=XA^{+\infty}$, and $\la a^{-\infty}\ra=XA^{-\infty}$.
Consequently, we have
$(a^{+\infty}|a^{-\infty})=|X|\stackrel{\rm{def}}{=}|\shl\la a\ra|$, and $$
\rho(a^{+\infty},a^{-\infty})\stackrel{\rm{def}}{=}\frac{1}{(a^{+\infty}|a^{-\infty})+1}=\frac{1}{|\shl\la
a\ra|+1}.$$ Now, we observe that $|a|=2|X|+|A|$, whence we get
$$(a|a^{+\infty})=(XAX^{-1}|XA^{+\infty})=|XA|=|X|+|A|>
\frac{2|X|+|A|}{2}=\frac{|a|}{2}.$$ Therefore, we have
$$\rho(a,a^{+\infty})\stackrel{\rm{def}}{=}\frac{1}{(a|a^{+\infty})+1}<
\frac{2}{|a|}.$$
\end{proof}

\section{The free group. Lemmas}

In this section, we prove several properties of the action of the free
non-abelian group $F$ on its boundary $\dd F$. We assume that $F$ has a fixed
system of generators, so we use the notation introduced in the previous
section.

\begin{proclaim}{\mmLF.1. Lemma}
Let $a\in{F}$, $k\in\N$. Assume that $|a|\ge 2k$. Then we have
$$aB_{1/k}(a^{-\infty})\ \cup\ B_{1/k}(a^{+\infty}) = \dd F.$$
\end{proclaim}

\begin{proof}
It suffices to show that for an arbitrary point $w\in\dd F$ we have either
$w\in B_{1/k}(a^{+\infty})$ or $a^{-1}w\in B_{1/k}(a^{-\infty})$.

The assumption $|a|\ge 2k$ and inequality~(\mmF.2) imply the following two
inequalities: $$(a^{+\infty}|a) > \frac{|a|}{2}\ge k; \eqno(\mmLF.1)\qquad $$ $$
(a^{-\infty}|a^{-1}) > \frac{|a^{-1}|}{2}=\frac{|a|}{2}\ge k. \eqno(\mmLF.2)\qquad
$$ Obviously, we can take an element $w'\in F$ so close to $w$ that the
equalities $(a^{-1}w|a^{-1})=(a^{-1}w'|a^{-1})$ and $(w|a)=(w'|a)$ hold. Then
we have
\begin{equation}\notag
\begin{split}
(a^{-1}w|a^{-1})+(w|a)&=(a^{-1}w'|a^{-1})+(w'|a)\\
&=\frac12\big(|a^{-1}w'|+|a^{-1}|-|w'|\big)+\frac12\big(|w'|+|a|-|a^{-1}w'|\big)\\
&=\frac12\big(|a^{-1}|+|a|\big)=|a|\ge 2k.
\end{split}
\end{equation}
This means that either
$(a^{-1}w|a^{-1}) \ge k$ or $(w|a) \ge k$. We observe that for every $x,y,z\in
F\cup\dd F$ the conditions $(x|y)\ge k$ and $(x|z)\ge k$ imply the condition
$(y|z)\ge k$ (see the properties of the Gromov product in Subsection~\mmF.4).
Consequently, if $(a^{-1}w|a^{-1}) \ge k$, then by~(\mmLF.2) we get $(a^{-1}
w|a^{-\infty})\ge k$, i.e., $a^{-1}w\in B_{1/(k+1)}(a^{-\infty})\subset
B_{1/k}(a^{-\infty})$. In the case where $(w|a) \ge k$, the inequality~(\mmLF.1)
yields $(w|a^{+\infty})\ge k$ whence $w\in B_{1/k}(a^{+\infty})$, as required.
\end{proof}

\begin{proclaim}{\mmLF.2. Lemma}
Let $a\in{F}$, and let $|\shl(a)|\ge k\in\N$. Then the $k$-element collection
$\{a,a^2,\dots,a^k\}\subset F$ is totally $1/k$-contracting for the metric
$F$-space $(\dd F,\rho)$.
\end{proclaim}

\begin{proof}
It is necessary (and it is sufficient) to prove that for every measure
$\lambda\in\PP(\dd F)$ there are a point $w\in\dd F$ and an element $g\in\{a,
a^2, \ldots, a^k\}$ such that $$ g\lambda\big(B_{1/k}(w)\big)\ge 1-{1/k}. $$
The assumption $|\shl(a)|\ge k$ and~(\mmF.1) imply that
$\rho(a^{-\infty},a^{+\infty})< 1/k$, whence it follows that
$B_{1/k}(a^{-\infty})=B_{1/k}(a^{+\infty})$ (see Subsection~\mmF.4). We denote
the ball $B_{1/k}(a^{-\infty})=B_{1/k}(a^{+\infty})$ by $B_{1/k}$. Since $|a|>
2|\shl(a)|\ge 2k$, it follows by Lemma~\mmLF.1 that for every $i\neq j\in\Z$
we have $$ a^{|i-j|} B_{1/k}\cup B_{1/k}=\dd F, $$ which is equivalent to $$
a^{i} B_{1/k} \cup a^{j} B_{1/k}=\dd F. $$

Thus, the sets in the collection $\{\dd F\ssm a^z B_{1/k}\}_{z\in\Z}$ are
pairwise disjoint. Hence the set $\{1,\ldots,k\}$ contains a number ($i_0$
say) such that $$ \lambda(\dd F\ssm a^{-i_0} B_{1/k})\le 1/k. $$ Also, we
observe that $$ \lambda(\dd F\ssm a^{-i_0} B_{1/k})=1-\lambda(a^{-i_0}
B_{1/k})=1-a^{i_0}\lambda(B_{1/k}). $$ Combining the latter two formulas we
get $$ a^{i_0}\lambda(B_{1/k}) \ge 1-1/k. $$ Thus, for the element
$g:=a^{i_0}\in\{a, a^2, \ldots, a^k\}$ and for the point $a^{+\infty}$ we have
$$g\lambda\big(B_{1/k}(a^{+\infty})\big)\ge 1-{1/k},$$ as required.
\end{proof}

\begin{proclaim}{\mmLF.3. Lemma (cf.~\cite[Proposition~2.27]{LSh})}
Assume that elements $a$ and $b$ of the free
group $F$ do not commute. Let $k\in\N$. Then the set $\{a,\ b,\ a^{20k}ba^{-20k},\
b^{20k}ab^{-20k}\}$ has an element $h$ with $|\shl\la h\ra|\ge k$.
\end{proclaim}

\begin{proof}
Recall that $\la x\ra$ denotes the reduced word representing an element $x\in
F$. For a word $W\in U^*$ we will denote by $\wh{W}$ the element of $F$ represented by
$W$. (We thus have $\wh{\la x\ra}=x$.)

Suppose that $\la a\ra=XAX^{-1}$ and $\la b\ra=YBY^{-1}$, where $X$, $Y$ are
reduced words, and $A$, $B$ are cyclically reduced words (as mentioned in
Subsection~\mmF.1\ above, such $X$, $Y$, $A$ and $B$ exist and unique).
Without loss of generality, we assume that $|A|\le|B|$. We will show that
under this assumption at least one of the elements $a$, $b$, and
$b^{20k}ab^{-20k}$ has the wing of length $\ge k$.

Assume the converse, i.e., assume that we have $$ |X|<k, \qquad |Y|<k, \qquad
|\shl\la b^{20k}ab^{-20k}\ra|<k. $$ We observe that $$
b^{20k}ab^{-20k}=\wh{Y}\wh{B}^{20k}\wh{Y}^{-1}\wh{X}\wh{A}\wh{X}^{-1}\wh{Y}\wh{B}^{-20k}\wh{Y}^{-1}.
$$ Let $\qq$ denote the element $\wh{Y}^{-1}\wh{X}\wh{A}\wh{X}^{-1}\wh{Y}$.
Then $$ b^{20k}ab^{-20k}=\wh{Y}\wh{B}^{20k}\qq \wh{B}^{-20k}\wh{Y}^{-1}. $$
Now, we concentrate on properties of the element $\wh{B}^{20k}\qq$ and of the
corresponding reduced word $\la\wh{B}^{20k}\qq\ra$. Let us prove several
inequalities.

Since $\qq=\wh{Y}^{-1}\wh{X}\wh{A}\wh{X}^{-1}\wh{Y}$, the inequalities
$|Y|<k\le k|B|$, $|X|<k\le k|B|$, and $|A|\le |B|\le k|B|$ imply the following
inequality: $$ |\qq|<5k|B|. \eqno(\mmLF.3)\qquad $$ Then, since
$|\wh{B}^{20k}|=20k|B|$, we have $$ |\wh{B}^{20k}\qq|<25k|B|. \eqno(\mmLF.4)\qquad
$$ Since the lengths of cores of conjugate elements are equal, it follows that
$|\osn \la b^{20k}ab^{-20k} \ra|=|A|$. By assumption, we have $|\shl\la
b^{20k}ab^{-20k}\ra|<k$, whence it follows that $$ |\wh{Y}\wh{B}^{20k}\qq
\wh{B}^{-20k}\wh{Y}^{-1}|=|b^{20k}ab^{-20k}|<|A|+2k. $$ Then, since $|Y|<k$,
we have $$ |\wh{B}^{20k}\qq \wh{B}^{-20k}| < |A|+2k+2k\le 5k|B|.
\eqno(\mmLF.5)\qquad $$ Now, we estimate the Gromov product of the elements
$(\wh{B}^{20k}\qq)^{-1}$ and $\wh{B}^{-20k}$. We observe that
\begin{equation}\notag
\begin{split}
\Big(\big(\wh{B}^{20k}\qq\big)^{-1}\Big|\wh{B}^{-20k}\Big)&=\frac12\Big(|\wh{B}^{20k}\qq|+|\wh{B}^{20k}|-|\wh{B}^{20k}\qq\wh{B}^{-20k}|\Big)\\
&\ge \frac12\Big(|\wh{B}^{20k}|-|\qq|+|\wh{B}^{20k}|-|\wh{B}^{20k}\qq\wh{B}^{-20k}|\Big)\\
&=20k|B|-\frac12\Big(|\qq|+|\wh{B}^{20k}\qq\wh{B}^{-20k}|\Big).
\end{split}
\end{equation}
Then by~(\mmLF.3) and by~(\mmLF.5) we obtain
$$\Big(\big(\wh{B}^{20k}\qq\big)^{-1}\Big|\wh{B}^{-20k}\Big)>15k|B|.$$
This means (see Subsection~\mmF.4) that the $15 k|B|$-prefixes of words
$\la\wh{B}^{20k}\qq\ra^{-1}$ and $B^{-20k}$ coincide, hence the word $B^{-15k}$ is
a prefix of $\la \wh{B}^{20k}\qq\ra^{-1}$. Equivalently, the word $B^{15k}$ is
a suffix of $\la \wh{B}^{20k}\qq\ra$.

At the other hand, since the word $B^{20k}$ of length $20k|B|$ is reduced and
since $|\qq|<5k|B|$ (see (\mmLF.3)), it follows that $B^{15k}$ is a prefix of
$\la\wh{B}^{20k}\qq\ra$.

We have thus shown that the word $\la \wh{B}^{20k}\qq\ra$ has suffix $B^{15k}$
and prefix $B^{15k}$, and, moreover, $|\la\wh{B}^{20k}\qq\ra|\le 25k|B|$
(see~(\mmLF.4)). Consequently, since the word $B$ is cyclically reduced, it
follows that each of the words $B \la B^{20k}\qq\ra$ and $\la B^{20k}\qq\ra
B$:

\begin{itemize}

\item[--] is reduced;

\item[--] has prefix $B^{15k}$ of length $15k|B|$;

\item[--] has suffix $B^{15k}$ of length $15k|B|$;

\item[--] has length at most $25k|B|+|B|\le 26k|B|$.

\end{itemize}

Therefore (since $|B^{15k}|+|B^{15k}| > 26k|B|$), the words $B \la
B^{20k}\qq\ra$ and $\la B^{20k}\qq\ra B$ coincide. Hence the elements
$\wh{B}^{20k}\qq$ and $\wh{B}$ of $F$ commute. Then the elements $\qq$ and
$\wh{B}$ also commute, whence $\wh{Y}\qq \wh{Y}^{-1}$ and
$\wh{Y}\wh{B}\wh{Y}^{-1}$ commute. We observe that $\wh{Y}\qq
\wh{Y}^{-1}=\wh{X}\wh{A}\wh{X}^{-1}=a$ and $\wh{Y}\wh{B}\wh{Y}^{-1}=b$. We
have thus proved that $a$ and $b$ commute, which contradicts the assumption of
our lemma.
\end{proof}

\begin{proclaim}{\mmLF.4. Corollary}
Assume that elements $a$ and $b$ of the free
group $F$ do not commute. Let $k\in\N$. Then the collection $Q^{50k^2}$, where $$ Q:=\{a,\
b,\ a^{-1},\ b^{-1},\ e\}, $$ is totally $1/k$-contracting for the metric
$F$-space $(\dd F,\rho)$.
\end{proclaim}

\begin{proof}
Since $Q$ contains the group identity, it follows that for every $r,s\in \N$
with $r\le s$ we have $Q^r\subset Q^s$. This implies that the elements $a$,
$b$, $a^{20k}ba^{-20k}$, and $\ b^{20k}ab^{-20k}$ lie in the set $Q^{40k+1}$,
whence by Lemma~\mmLF.3 it follows that the set $Q^{40k+1}$ contains an
element (say, $h$) with $|\shl\la h\ra|\ge k$. Then by Lemma~\mmLF.2 the
collection $\{h, h^2, \ldots, h^k\}$ is totally $1/k$-contracting for $(\dd
F,\rho)$. It is clear that $$ \{h, h^2, \ldots, h^k\}\subset
\big(Q^{40k+1}\big)^k\subset Q^{50k^2}. $$ We have thus proved that the
collection $Q^{50k^2}$ contains a totally $1/k$-contracting subcollection.
Hence, $Q^{50k^2}$ is totally $1/k$-contracting.
\end{proof}

\section{Proof of Theorem~1}

The outline of the proof is as follows: using the results of Section~\mmLF, we
show that the action of $G$ on the boundary $\dd F$ satisfies the sufficient
condition of $\mu$-proximality, which was described in Section~\mmPRIZ.

We assume that $F$ has a fixed system of generators, so we enable the notation
introduced in Sections~\mmF\ and~\mmLF.

\begin{proclaim}{Claim}
For an arbitrary large $k\in\N$, the group $G$ contains a finite universally
$1/k$-contracting collection for the metric $G$-space $(\dd F,\rho)$.
\end{proclaim}

\begin{proof}[Proof of the claim]
Let $x,y\in F\subset G$ be a pair of non-commuting elements. We will show that
the finite collection $$ R:=\{x,\ y,\ x^{-1},\ y^{-1},\ e\}^{50k^2}\subset
F\subset G $$ is universally $1/k$-contracting for $(\dd F,\rho)$. In view of
Remark~\mmPRIZ.3, it suffices to show that for every $g\in G$ the collection
$gRg^{-1}$ is totally $1/k$-contracting. We put $a:=gxg^{-1}$ and
$b:=gyg^{-1}$. Since $F$ is a normal subgroup of $G$, it follows that $a$ and
$b$ lie in $F$. We also observe that $a$ and $b$ do not commute (because $x$
and $y$ do not commute). It is clear that $$ gRg^{-1}=\{a,\ b,\ a^{-1},\
b^{-1},\ e\}^{50k^2}. $$ Thus, the collection $gRg^{-1}$ is totally
$1/k$-contracting by Corollary~\mmLF.4. The claim is proved.
\end{proof}

By the above claim, the action of $G$ on $\dd F$ meets the assumptions of
Proposition~\mmPRIZ.4. Therefore, the boundary $\dd F$ is a $\mu$-proximal
space. This means that there exists a unique $\mu$-stationary measure $\nu$ on
$\dd F$ (see Lemma~\mmRWL.3), and the pair $(\dd F,\nu)$ is a $\mu$-boundary of
$(G,\mu)$. To complete the proof, we observe that for each point $w\in\dd F$
the orbit $G(w)\supset F(w)$ is infinite, so $\nu$ is continuous by
Lemma~\mmRWL.4. Theorem~1 is proved.

\section{A theorem about the selective convergence}

\paragraph{Definition \rm{(the selective convergence)}.}
Let $T$ be a topological space. We say that a sequence
$\{(a_i,b_i)\}_{i\in\Z_+}$, where $(a_i,b_i)$ is a pair of points
in $T$, {\it selectively converges\/} to a point $t\in T$ if there
exist a converging to $t$ sequence $\{c_i\}_{i\in\Z_+}$ such that
$c_i\in \{a_i,b_i\}$ for each~$i$. Obviously, a sequence of pairs
of points in a Hausdorff space can selectively converge to at most
two distinct points.

\medskip

In this section, we prove the following theorem.

\begin{proclaim}{Theorem about the selective convergence}
Let $G$ be a countable group with a normal free non-abelian subgroup $F$. Let
$\mu$ be a nondegenerate measure on $G$. Then for a.e. path
$\tau=\{\tau_i\}_{i\in\Z_+}$ of the random $\mu$-walk and for an arbitrary
nontrivial element $a\in F$, the sequence of pairs $\{(\tau_i a \tau_i^{-1},
\tau_i a^{-1}\tau_i^{-1})\}_{i\in\Z_+}$ selectively converges to a point
$w(\tau)\in\dd F$.
\end{proclaim}

\begin{proof}
We assume that $F$ has a fixed system of generators, and we consider the
corresponding word metric, Gromov product, and metric $\rho$ in $F\cup\dd F$.
As above, for $v\in\dd F$ and $\varepsilon>0$, we set
$B_\varepsilon(v):=\{w\in \dd F|\ \rho(v,w)\le\varepsilon\}$ (see the notation
in Section~\mmF).

The image of a point $v\in \dd F$ under the action of an element $g\in G$ we
denote by $gv$. We note that for every nontrivial element $x\in F$ and every
$g\in G$ the point $gx^{+\infty}$ coincides with the point
$(gxg^{-1})^{+\infty}$.

It is clear that under the assumptions of the theorem our Theorem~1 is
applicable, thus on $\dd F$ there exists a unique $\mu$-stationary measure
$\nu$.

\begin{proclaim}{Claim~(1)}
Let $k\in\N$, and $b\in F$ be a nontrivial element. Then
there exists a constant $C(b,k)>0$ such that for each $g\in G$ we have $$
g\nu\, \Big(\,B_{1/k}(g b^{+\infty})\ \cup\ B_{1/k}(g b^{-\infty})\,\Big) \ge
C(b,k). $$
\end{proclaim}

\noindent{\it Remark.\/} This claim is not trivial, because we may have
$B_{1/k}(g b^{\pm\infty})\neq g B_{1/k}(b^{\pm\infty})$.

\begin{proof}[Proof of Claim~$(1)$]
Let $g$ be an arbitrary element of $G$. We put $$ E:=B_{1/k}(g
b^{+\infty})\cup B_{1/k}(g b^{-\infty}). $$ Obviously, for every element $x\in
F$ and every positive integer $r$, the word $\osn\la x^r\ra$ coincides with
the word $(\osn\la x\ra)^r$. Therefore, since the element $gbg^{-1}\in F$ is
not trivial, we have $$ |(gbg^{-1})^{2k}|\ge 2k|\osn(gbg^{-1})| \ge 2k. $$ For
every $r\in\N$, we obviously have $(gb^{r}g^{-1})^{\pm\infty}=g
(b^r)^{\pm\infty}=g b^{\pm\infty}$, whence $B_{1/k}(g
b^{\pm\infty})=B_{1/k}\big((gb^{r}g^{-1})^{\pm\infty}\big)$. Then by
Lemma~\mmLF.1 we have $$ (gbg^{-1})^{2k}E \cup E = \dd F. $$ Consequently, $$
b^{2k}g^{-1}E \cup g^{-1}E = g^{-1}\dd F=\dd F. $$ Applying the second part of
Lemma~\mmRWL.5, we obtain $$ g\nu(E)=\nu(g^{-1}E)\ge C''_{b^{2k}}. $$ To
complete the proof, we set $C(b,k)=C''_{b^{2k}}$.
\end{proof}

\begin{proclaim}{Claim~(2)}
If a sequence of measures $\{g_i (\nu)\}_{i\in\Z_+}$, where $g_i\in G$,
converges to a point measure $\delta_v$ with $v\in\dd F$, then for each
nontrivial element $b\in F\ssm e$ the sequence of pairs $\{(g_i b^{+\infty},
g_i b^{-\infty})\}_{i\in\Z_+}$ selectively converges to the point~$v$.
\end{proclaim}

\begin{proof}[Proof of Claim~$(2)$]
It is sufficient to show that for each $k\in\N$ there exists $N_{1/k}\in\N$
such that for each $j\ge N_{1/k}$ one of the points of the pair  $(g_j
b^{+\infty}, g_j b^{-\infty})$ lies at distance at most $1/k$ from~$v$ (with
respect to the metric $\rho$). Let us prove it. Let $k\in\N$. By Claim~(1)
there exists $C(b,k)>0$ such that for each $g\in G$ we have $$ g\nu\,
\Big(\,B_{1/k}(g b^{+\infty})\ \cup\ B_{1/k}(g b^{-\infty})\,\Big) \ge C(b,k).
$$ Since the sequence of measures $\{g_i (\nu)\}_{i\in\Z_+}$ converges to a
point measure $\delta_v$, it follows that there exists $N_{1/k}\in\N$ such
that for each $j\ge N_{1/k}$ we have $$ g_j\nu\Big(B_{1/k}(v)\Big)> 1-C(b,k).
$$ Thus we have $$ g_j\nu\, \Big(\,B_{1/k}(g_j b^{+\infty})\ \cup\ B_{1/k}(g_j
b^{-\infty})\,\Big) + g_j\nu\Big(B_{1/k}(v)\Big) > 1. $$ Consequently, the
ball $B_{1/k}(v)$ must intersect at least one of the balls
$B_{1/k}(g_jb^{+\infty})$ and $B_{1/k}(g_j b^{-\infty})$. This means, as
observed in Section~\mmF, that at least one of the points $g_jb^{+\infty}$,
$g_jb^{-\infty}$ lies at distance at most $1/k$ from the point~$v$. So, the
claim is proved.
\end{proof}

We denote by $\dd_rF$ the set of ``{\it rational\/}'' points of the boundary
$\dd F$: $$ \dd_rF:=\{x^{+\infty}|\ x\in F\ssm e\}\subset\dd F. $$

\begin{proclaim}{Claim~(3)}
If a sequence of measures $\{g_i (\nu)\}_{i\in\Z_+}$, where $g_i\in G$,
converges to a point measure $\delta_v$ with $v\in\dd{F}\ssm\dd_r{F}$, then for each
nontrivial element $b\in F\ssm e$ the sequence of numbers $\{|g_i b
g_i^{-1}|\}_{i\in\Z_+}$ tends to infinity.
\end{proclaim}

\begin{proof}[Proof of Claim~$(3)$]
Assume converse. Then the sequence $\{g_i b g_i^{-1}\}_{i\in\Z_+}$ has a
constant subsequence $\{g_{i_k} b g_{i_k}^{-1}\}_{k\in\Z_+}=y,y,y,\ldots$. By
Claim~(2), the sequence of pairs $\{(g_i b^{+\infty}, g_i
b^{-\infty})\}_{i\in\Z_+}$ selectively converges to the point~$v$, whence the
subsequence $\{(g_{i_k} b^{+\infty}, g_{i_k} b^{-\infty})\}_{k\in\Z_+}$ $=$
$\{(y^{+\infty}, y^{-\infty})\}_{k\in\Z_+}$ also selectively converges to~$v$.
But this means that $v\in\{y^{+\infty},y^{-\infty}\}\subset\dd_rF$, which
contradicts the assumption that $v\in\dd F\ssm\dd_rF$. Claim~(3) is thus
proved.
\end{proof}

\begin{proclaim}{Claim~(4)}
If a sequence of measures $\{g_i (\nu)\}_{i\in\Z_+}$, where $g_i\in G$,
converges to a point measure $\delta_v$ with $v\in\dd F\ssm \dd_rF$, then for
each nontrivial element $b\in F\ssm e$ the sequence of pairs $\{(g_i b
g_i^{-1}, g_i b^{-1}g_i^{-1})\}_{i\in\Z_+}$ selectively converges to~$v$.
\end{proclaim}

\begin{proof}[Proof of Claim~$(4)$]
From Claim~(2) it follows that there exists a sequence
$\{\delta_i\}_{i\in\Z_+}$, where $\delta_i\in\{+1,-1\}$, such that the
sequence $\{g_i (b^{\delta_i})^{+\infty}\}_{i\in\Z_+}$ converges to~$v$. By~(\mmF.3),
for each $i$ we have $$ \rho(g_i b^{\delta_i}
g_i^{-1}, g_i (b^{\delta_i})^{+\infty})< \frac{2}{|g_i b^{\delta_i}
g_i^{-1}|}. $$ By Claim~(3) we have $|g_i b^{\delta_i}g_i^{-1}|\to+\infty$,
whence $$ \rho (g_i b^{\delta_i} g_i^{-1}, g_i (b^{\delta_i})^{+\infty})\to 0.
$$ This implies that the sequence $\{g_i b^{\delta_i} g_i^{-1}\}_{i\in\Z_+}$
converges to~$v$ (because the sequence $\{g_i
(b^{\delta_i})^{+\infty}\}_{i\in\Z_+}$ converges to~$v$), which means by
definition that the sequence of pairs $\{(g_i b g_i^{-1}, g_i
b^{-1}g_i^{-1})\}_{i\in\Z_+}$, indeed, selectively converges to~$v$.
\end{proof}

To complete the proof of our theorem, we observe that by Theorem~1 the
$\mu$-stationary measure $\nu$ is continuous and the pair $(\dd F,\nu)$ is a
$\mu$-boundary of $(G,\mu)$. The latter fact means that for a.e. path
$\tau=\{\tau_i\}_{i\in\Z_+}$ of the random $\mu$-walk the sequence of measures
$\{\tau_i (\nu)\}_{i\in\Z_+}$ converges to a point measure $\delta_{w(\tau)}$,
where $w(\tau)\in\dd F$. Furthermore, since the set $\dd_rF$ is countable and
since the measure $\nu$ is continuous, it follows that $\nu(\dd_rF)=0$. From
this we conclude by Theorem~\mmRWL.2 that {\it for a.e. path $\tau$ we have}
$w(\tau)\in\dd F\ssm\dd_rF$. Therefore, the theorem is proved now by
Claim~(4).
\end{proof}

\section{Proof of Theorem~2}

By the assumption of Theorem~2, $F$ contains
an element $\u\neq e$ that is fixed under the action of the subgroup $A$.
By the selective convergence theorem (see Section~\mmVYB), for a.e. path
$\tau=\{\tau_i\}_{i\in\Z_+}=\{x_i\alpha_i\}_{i\in\Z_+}$ the sequence of pairs
$\{(\tau_i \u \tau_i^{-1}, \tau_i \u^{-1}\tau_i^{-1})\}_{i\in\Z_+}$
selectively converges to a point~$w(\tau)\in\dd F$. Since $\u$ commutes with
the elements of~$A$, it follows that $$ \tau_i \u \tau_i^{-1}=x_i \u x_i^{-1}
\quad\text{and}\quad \tau_i \u^{-1} \tau_i^{-1}=x_i \u^{-1} x_i^{-1}. $$
Therefore, the sequence of pairs $\{(x_i \u x_i^{-1}, x_i
\u^{-1}x_i^{-1})\}_{i\in\Z_+}$ selectively converges to~$w(\tau)$. This means
that there exists a sequence $\{\delta_i\}_{i\in\Z_+}$ with
$\delta_i\in\{+1,-1\}$ such that the sequence $\{x_i \u^{\delta_i}
x_i^{-1}\}_{i\in\Z_+}$ converges to~$w(\tau)$. In particular, if we fix in $F$
a system of generators and consider the associated word metric, then
$|x_i \u^{\delta_i} x_i^{-1}|\to\infty$ as $i\to\infty$. Also, we observe that
for the Gromov products $(x_i|x_i \u^{\delta_i} x_i^{-1})$ the following
inequality holds: $$ (x_i|x_i \u^{\delta_i} x_i^{-1})=\frac12\big(|x_i|+|x_i
\u^{\delta_i} x_i^{-1}|-|x_i \u^{-\delta_i}|\big)\ge \frac12\big(|x_i
\u^{\delta_i} x_i^{-1}|-|\u|\big). $$ Thus, the sequence $\{(x_i|x_i
\u^{\delta_i} x_i^{-1})\}_{i\in\Z_+}$ tends to infinity. Then, since the
sequence $\{x_i \u^{\delta_i} x_i^{-1}\}_{i\in\Z_+}$ converges to~$w(\tau)$, it
follows that the sequence $\{x_i\}_{i\in\Z_+}$ also converges to~$w(\tau)$.
Theorem~2 is thus proved.

\section{The stability of the Markov--Ivanovsky normal form}

\begin{proof}[Proof of Theorem~{\rm3}]
In the proof, we use the notation introduced in Section~\mmBG. We recall that
the pure braid group $P_n$ is a semidirect product of the normal free subgroup
$F_{n-1}$ by the subgroup $P_{n-1}$, and each element $\gamma \in P_n$ can be
written as a unique product $\gamma = x \alpha$, where $x \in F_{n-1}$ and
$\alpha \in P_{n-1}$. It is easy to check that the element $$
\u=s_{n1}s_{n2}\cdots
s_{n(n-1)}=\si_{n-1}\ldots\si_2\si_1\si_1\si_2\cdots\si_{n-1}\in F_{n-1} $$
commutes with all elements of the subgroup~$P_{n-1}$. We can thus apply
Theorem~2 to the case of the pure braid group $P_n=F_{n-1} \rtimes P_{n-1}$. By
Theorem~2, for each nondegenerate measure $\mu\in\PP(P_n)$ and for a.e. path of
the random $\mu$-walk $\tau=\{\tau_i\}_{i\in\Z_+}=\{x_i \alpha_i\}_{i\in\Z_+}$
(here, $x_i\in F_{n-1}$, $\alpha_i\in P_{n-1}$), the sequence
$\{x_i\}_{i\in\Z_+}$ converges to a point on the boundary~$\dd F_{n-1}$. By the
definition of the Markov--Ivanovsky normal form, the reduced representative of
the element $x_i$ (we mean the reduced representative over the generators
$\{s_{nj}, 1\le j <n\}$ and their inverses) is a prefix of the normal form
$\mathfrak{I}_P(\tau_i)$. This means exactly that the normal form
$\mathfrak{I}_P$ is stable in the pure braid group $P_n$.

Since $P_n$ is a subgroup of finite index in $B_n$, the following lemma proves
the stability of the Markov--Ivanovsky normal form $\mathfrak{I}_B$ in the
braid group $B_n$.
\end{proof}

\begin{proclaim}{\nJ.1. Lemma\footnote{This lemma is an interpretation
of the well-known fact that the exit boundaries of a group and of its       
subgroup of finite index coincide. But due to the obvious technical reasons we
can not regard the lemma as just a corollary of this fact.}} Let $G$ be a
countable group, $Q$ be a subgroup of finite index in $G$, and $\Pi
\subset G$ be a set of left coset representatives for $Q$. Let
$\mathfrak{N}_Q$ be a normal form in the group $Q$. Denote by $\mathfrak{N}_G$
a normal form in $G$ that takes an element $g\in G$ to the word
$\mathfrak{N}_Q(g\pi_g^{-1})\pi_g$, where $\pi_g:=\Pi\cap Qg$. Suppose that
$\mathfrak{N}_Q$ is stable. Then $\mathfrak{N}_G$ is stable.
\end{proclaim}

\begin{proof}
By the definition of stable normal forms, it is necessary and sufficient to
prove that $\mathfrak{N}_G$ is $\mu$-stable for an arbitrary nondegenerate
measure $\mu$ on $G$.

For a path $\tau=\{\tau_i\}_{i\in\Z_+}$ and a subset $H\subset G$, we will
denote by $\tau^{\scriptscriptstyle{H}}$ the subsequence in $\tau$ that
consists of elements of $H$. Lemma~\mmRWL.6 implies that if $H$ is cofinite,
then $\tau^{\scriptscriptstyle{H}}$ is infinite for a.e. path $\tau$. We
observe that for each $g\in G$ the subset $Qg$ is cofinite in $G$, because $Q$
has finite index in $G$. Consequently, for each $g\in G$ and for $\mi$-a.e.
path $\tau$ the subsequence $\tau^{\scriptscriptstyle{Qg}}$ is infinite.

First, we show that for $\mi$-a.e. path $\tau$ the sequence of words
$\mathfrak{N}_Q(\tau^{\scriptscriptstyle{Q}})$ converges (throughout the
proof, for a normal form $\mathbf{N}$ and a sequence
$\mathbf{s}:=\{s_j\}_{j\in\Z_+}$, we denote by $\mathbf{N}(\mathbf{s})$ the
sequence $\{\mathbf{N}(s_j)\}_{j\in\Z_+}$). Let us define a measure $\mu'$ on
the group $Q$ by setting $\mu'(q)$ equal to the probability that for a sample
path $\tau$ of the random $\mu$-walk on $G$ we have
$\{\tau^{\scriptscriptstyle{Q}}\}_1=q$. Then Lemma~\mmRWL.6 implies that
$\mu'$ is a probability measure on $Q$. Since $\mu$ is nondegenerate for $G$,
it easily follows that $\mu'$ is nondegenerate for $Q$. Now, we consider the
random $\mu'$-walk on $Q$. Let $P_{\mu'}$ denote the corresponding measure on
the path space $Q^{\Z_+}$. Clearly, the ($\mi$-a.e. defined) mapping
$\tau\mapsto\tau^{\scriptscriptstyle{Q}}$ from $G^{\Z_+}$ to $Q^{\Z_+}$ maps
the measure $\mi$ to the measure $P_{\mu'}$. Since $\mathfrak{N}_Q$ is stable
(in $Q$) and $\mu'$ is nondegenerate (in $Q$), it follows that for
$P_{\mu'}$-a.e. path $\kappa=\{\kappa_i\}_{i\in\Z_+}$ of the random
$\mu'$-walk the sequence of words $\mathfrak{N}_Q(\kappa_i)$ converges.
Consequently, for $\mi$-a.e. path $\tau$ the sequence
$\mathfrak{N}_Q(\tau^{\scriptscriptstyle{Q}})$ converges.

From the above it follows by standard arguments that for $\mi$-a.e. path
$\tau$ and for each element $g\in G$ the sequence of words
$\mathfrak{N}_Q(\tau^{\scriptscriptstyle{Qg}}g^{-1})$ converges.

Now, we observe that for every elements $g,h\in G$ and for $\mi$-a.e. path
$\tau$ the sequences $\mathfrak{N}_Q(\tau^{\scriptscriptstyle{Qg}}g^{-1})$ and
$\mathfrak{N}_Q(\tau^{\scriptscriptstyle{Qh}}h^{-1})$ converge to one and the
same limit. (This follows from Lemma~\mmRWL.6, which implies that for
$\mi$-a.e. path $\tau=\{\tau_i\}_{i\in\Z_+}$ there exist $k\in\N$ and an
infinite set $L\subset\N$ such that for each $l\in L$ we have $\tau_l\in Qg$
and $\tau_{l+k}=\tau_l g^{-1}h\in Qh$, whence it follows that the sequences
$\tau^{\scriptscriptstyle{Qg}}g^{-1}$ and $\tau^{\scriptscriptstyle{Qh}}h^{-1}$
have infinite coinciding subsequences.)

To complete the proof, we observe that each path $\tau$ in $G^{\Z_+}$ splits
into a finite number of subsequences of the form
$\tau^{\scriptscriptstyle{Qp}}$, where $p\in\Pi$. Accordingly, the sequence of
words $\{\mathfrak{N}_G(\tau_i)\}_{i\in\Z_+}$ (which by the definition of
$\mathfrak{N}_G$ equals the sequence
$\{\mathfrak{N}_Q(\tau_i\pi_{\tau_i}^{-1})\pi_{\tau_i}\}_{i\in\Z_+}$) splits
into the subsequences $\mathfrak{N}_Q(\tau^{\scriptscriptstyle{Qp}}p^{-1})p$,
$p\in\Pi$. By the above, for $\mi$-a.e. path $\tau$ these subsequences
converge. Furthermore, they converge to one and the same (infinite) limit.
This means that for $\mi$-a.e. path $\tau$ the sequence of words
$\{\mathfrak{N}_G(\tau_i)\}_{i\in\Z_+}$ converges, i.e., the form
$\mathfrak{N}_G$ is $\mu$-stable, as required.
\end{proof}


\end{document}